\documentclass{amsart} 
\usepackage{fancyhdr}
\usepackage{amsmath,amsthm,amscd,amsfonts,amssymb}

\input {xy}

\xyoption{all}
\usepackage[all]{xypic}
\usepackage{latexsym}
\usepackage{color} 
\theoremstyle{plain}

\begin{document}
\input{amssym.def}

\numberwithin{equation}{section}

\newtheorem{guess}{Theorem}[section]
\newcommand{\bth}{\begin{guess}$\!\!\!${\bf }~}
\newcommand{\eeth}{\end{guess}}

\newtheorem{propo}[guess]{Proposition}

\newcommand{\bprop}{\begin{propo}$\!\!\!${\bf }~}
\newcommand{\eprop}{\end{propo}}

\newtheorem{lema}[guess]{Lemma}
\newcommand{\blem}{\begin{lema}$\!\!\!${\bf }~}
\newcommand{\elem}{\end{lema}}

\newtheorem{defe}[guess]{Definition}
\newcommand{\bdefe}{\begin{defe}$\!\!\!${\bf }~}
\newcommand{\edefe}{\end{defe}}

\newtheorem{coro}[guess]{Corollary}
\newcommand{\bcor}{\begin{coro}$\!\!\!${\bf }~}
\newcommand{\ecor}{\end{coro}}

\newtheorem{rema}[guess]{\it Remark}
\newcommand{\brem}{\begin{rema}$\!\!\!${\it }~\rm}
\newcommand{\erem}{\end{rema}}

\theoremstyle{remark}
\newtheorem{note}[guess]{Notation}

\newcommand{\pf}{{\noindent \bf Proof: }}
\newcommand{\enpf}{\begin{flushright} {\it q.e.d} \end{flushright}}

\newtheorem{eq}[guess]{equation}
\newcommand{\beq}{\begin{equation}}
\newcommand{\eeq}{\end{equation}}

\newtheorem{eqa}[guess]{eqnarray}
\newcommand{\beqa}{\begin{eqnarray}}
\newcommand{\eeqa}{\end{eqnarray}}

\newcommand{\G}{\Gamma}
\newcommand{\hra}{\hookrightarrow}
\newcommand{\lr}{\longrightarrow}
\newcommand{\K}{\mathrm K}
\newcommand{\mk}{\mathrm k}
\newcommand{\Ha}{\mathcal H}
\newcommand{\Rss}{\mathcal R^\Gamma}
\newcommand{\Gt}{\tilde \mathcal F ^ \Gamma}
\newcommand{\U}{\mathcal U}
\newcommand{\N}{\mathcal N}
\newcommand{\Ga}{\mathcal G}
\newcommand{\m}{\mathcal}
\newcommand{\g}{\goth }

\newcommand{\Vls}{\mathcal V \mid_{q_{\underline l}^{-1}(\underline s)}}
\newcommand{\Vms}{\mathcal V \mid_{q_{\underline m}^{-1}(\underline s)}}
\newcommand{\Vlas}{\mathcal V \mid_{q_{(\underline l)}^{-1}(\underline s)}}
\newcommand{\Vmas}{\mathcal V \mid_{q_{((\underline m))}^{-1}(\underline s)}}
\newcommand{\Vdms}{\mathcal V \mid_{q_{\underline m}^{-1}(\underline s)}^{**}}
\newcommand{\Vym}{\mathcal V \mid_{Y_{\underline m}}}
\newcommand{\Vyl}{\mathcal V \mid_{Y_{\underline l}}}

\newcommand{\Sma}{S_{m}^\Gamma}
\newcommand{\Sm}{S_{\underline m}^\Gamma}
\newcommand{\Sl}{S_{\underline l}^\Gamma}
\newcommand{\Sal}{S_{(\underline l)}^\Gamma}
\newcommand{\Sam}{S_{(\underline m)}^\Gamma}
\newcommand{\Smi}{S_{\underline m_i}^{\Gamma}}

\newcommand{\Zm}{Z_{\underline m}^{\Gamma}}

\newcommand{\Ym}{Y_{\underline m}^{\Gamma}}
\newcommand{\Yn}{Y_{\underline {m_0}}^{\Gamma}}
\newcommand{\Yam}{Y_{\underline (m)}^{\Gamma}}

\newcommand{\Yal}{Y_{\underline (l)}^{\Gamma}}

\newcommand{\Km}{K_{\underline m}^\Gamma}

\newcommand{\un}{\underline}
\newcommand{\us}{\underline s}

\newcommand{\qu}{q_{\underline m}}
\newcommand{\qm}{q_{{\underline m}_0}}

\newcommand{\pu}{p_{\underline m}}
\newcommand{\phim}{\phi_{\underline m}}

\newcommand{\gb}{$\Gamma$-bundle}
\newcommand{\gsb}{$\Gamma$-subbundle}

\newcommand{\qis}{q_{\underline m}^{-1}(\underline s)}

\newcommand{\Ul}{U_{\underline l}}
\newcommand{\Um}{U_{\underline m}}
\newcommand{\Uam}{U_{(\underline m)}}
\newcommand{\Ual}{U_{(\underline l)}}

\newcommand{\Wm}{\mathcal W_{\underline m}}

\newcommand{\Gm}{\tilde{\mathcal G_m}}
\newcommand{\Gn}{\tilde{\mathcal G_{m_0}}}
\newcommand{\pmV}{\phi_{\underline m}^*(p_{\underline m}^* V)}

\newcommand{\wpmV}{\mathbb P (\wedge^{r}\phi_{\underline
    m}^*(p_{\underline m}^* V))}

\newcommand{\pmwV}{\mathbb P (\phi_{\underline m}^*(p_{\underline
    m}^*(\wedge^{r}V)))}

\newcommand{\pV}{\phi_{\underline {m_0}}^*(p_{\underline {m_0}}^* V)}

\newcommand{\E}{\mathcal E}
\newcommand{\Ep}{\mathcal E'}
\newcommand{\Epp}{\mathcal E''}

\newcommand{\cu}{{\mathcal U}}
\newcommand{\cf}{{\mathcal F}}
\newcommand{\ce}{{\mathcal E}}
\newcommand{\co}{{\mathcal O}}
\newcommand{\cg}{{\mathcal G}}
\newcommand{\cm}{{\mathcal M}}

\newcommand{\mcp}{{\mathcal P}}
\newcommand{\cp}{{\sf P}}

\newcommand{\cq}{{\sf Q}}
\newcommand{\cs}{{\mathcal S}}
\newcommand{\cl}{{\mathcal L}}

\newcommand{\ctext}[1]{\makebox(0,0){#1}}
\setlength{\unitlength}{0.1mm}

\def\bC{{\Bbb C}}
\def\a{\alpha }
\def\bR{{\Bbb R}}
\def\tD{{\tilde D}}

\newcommand{\Mo}{M_\Gamma^{\mu ss}}
\newcommand{\Ro}{\R_\Gamma^{\mu ss}}
\newcommand{\Q}{Quot(\mathcal H, P)}
\newcommand{\gF}{gr^\mu_\G(F)}
\newcommand{\Fs}{gr^\mu_\G(F)^*}

\newcommand{\Di}{\D_{k_i}}
\newcommand{\DK}{\D_K}

\title [Donaldson-Uhlenbeck compactification]{Parabolic bundles on
  algebraic surfaces I- the Donaldson--Uhlenbeck compactification}

\author{V. Balaji}

\address{Chennai Mathematical Institute, 92, G.N. Chetty Road,
Chennai 600017, India}

\email{balaji@cmi.ac.in}

\author[A. Dey]{A. Dey}

\address{The Institute of Mathematical Sciences, CIT
Campus, Taramani, Chennai 600113, India}

\email{arijit@imsc.res.in}

\author{R. Parthasarathi}\thanks{The third author was supported by
    the National Board for Higher Mathematics, India}

\address{Chennai Mathematical Institute, 92, G.N. Chetty Road,
Chennai 600017, India}

\email{partha@cmi.ac.in}


\begin{abstract}
  The aim of this paper is to construct the parabolic version of the
  Donaldson--Uhlenbeck compactification for the moduli space of
  parabolic stable bundles on an algenraic surface with parabolic
  structures along a divisor with normal crossing singularities. We
  prove the non--emptiness of the moduli space of parabolic stable
  bundles of rank $2$ and also prove the existence of components with
  smooth points.
\end{abstract}

\maketitle

\section{Introduction}
Let $X$ be a smooth projective variety defined over the field $\bC$ of
complex numbers. Moduli spaces of sheaves with parabolic structures
were defined and constructed in great generality by Maruyama and
Yokogawa (\cite{MY}). This work of theirs generalises the earlier
construction of Mehta and Seshadri (\cite{MS}) when $dim(X) = 1$. When
$dim(X) = 2$, i.e $X$ is a smooth projective surface and if $D$ is an
effective divisor on $X$ then one finds from the work of Kronheimer
and Mrowka (cf \cite{kron} and \cite{kron1}) that the underlying geometry
and topology of the moduli space of parabolic bundles of rank two and
trivial determinant have very interesting applications arising out of
a generalization of Donaldson polynomials defined from these moduli
spaces. These moduli spaces and their compactifications were studied
in the papers of Kronheimer and Mrowka but primarily from the
differential geometric standpoint. In particular, the
Kobayashi-Hitchin correspondence was conjectured in these papers and
this has since been proven by a number of people in growing order of
generality. (cf \cite{biquard}, \cite{linarasimhan}, \cite{steer}).

The purpose of this paper and its sequel (\cite{balasuman}) is to
initiate a comprehensive study of the geometry of the moduli space of
$\mu$--stable parabolic bundles of arbitrary rank on {\it smooth
  projective surfaces} with parabolic structures on an reduced divisor
$D$ with {\it normal crossing singularities}. More precisely, in this
paper we construct the analogue of the Donaldson-Uhlenbeck
compactification of the moduli space of $\mu$--stable parabolic
bundles of arbitrary rank and also prove the existence of $\mu$-stable
parabolic bundles when certain topological invariants are allowed to
be arbitrarily large. We also show the existence of components with
smooth points. We summarise our results in the following theorem. For
notations see \eqref{modulinotation}:

\bth 
\begin{enumerate}
\item There exists a natural compactification of the moduli space
  $M^{\pmb{\alpha}}_{k, {\mathbf j}, \bf r}(r, \mathcal P, \kappa)$ of
  $\mu$--stable parabolic bundles with fixed determinant $\mathcal P$
  and with fixed topological and parabolic datum. Furthermore, the
  compactification can be set--theoretically be described as follows:

  \beqa {\overline {M^{\pmb{\alpha}}_{k, {\mathbf j}, \bf r}(r,
      \mathcal P, \kappa)}} \subset \coprod_{l\ge 0}
      M^{\pmb{\alpha}-poly}_{k', {\mathbf j}', \bf r}(r, \mathcal P,
      \kappa - l) \times S^l(X).  \eeqa
      
      where by, $M^{\pmb{\alpha}-poly}_{k, {\mathbf j}, \bf r}(r,
      \mathcal P, \kappa)$, we mean the set of isomorphism classes of
      {\it polystable} parabolic bundles with parabolic datum given by
      $(\pmb{\alpha},\bf l,\bf r, \bf j)$, fixed determinant $\mathcal
      P$ and with topological datum given by $k$ and $\kappa$.

    \item The moduli space of $\mu$--stable parabolic bundles of rank
      $2$ is non--empty, when the invariants $k$ and $\mathbf j$ are
      made sufficiently large and the weights satisfy some natural
      bounds. (see Theorem \ref{nonemptiness})
    \item Under these asymptotic assumptions, the moduli space has a
      component with smooth points.

\end{enumerate}
\eeth

This paper can therefore be seen as completing the algebro-geometric
analogue of the Kobayashi-Hitchin correspondence for parabolic bundles
on surfaces. We compare the moduli spaces that we construct with that
of Kronheimer-Mrowka when we restrict ourselves to the rank two case.

The main strategy used for the construction is to use the categorical
correspondence of the category of $\G$--bundles of fixed type $\tau$
on a certain Kawamata cover of the surface $X$ with the category of
parabolic bundles on $X$ with fixed parabolic datum (see \S1 for
definitions and terminology). The Kawamata cover $Y$ is non-canonical
and is therefore employed only as a stepping stone for the
construction.  Although non-canonical, the moduli problem gets defined
more naturally on $Y$ and one takes recourse to the ideas of Li and Le
Potier, as well as the earlier work of Donaldson to give an
algebraic--geometric construction of the Donaldson-Uhlenbeck
compactification of the moduli space of $\mu$--stable $\G$--bundles on
$Y$. Then by using the correspondence one can interpret the
compactification in a canonical manner as a compactification of the
moduli space of parabolic bundles over the surface $X$ with given
parabolic datum, thereby removing the non-canonical nature of the
construction. We believe that this moduli space can be realised, as in
the usual setting, as a generalized blow-down of the Maruyama-Yokogawa
moduli space. Unlike our moduli space, the Maruyama-Yokogawa space is
a GIT construction using Gieseker type stability for parabolic
sheaves.

We then go on to show that the moduli space of $\mu$--stable
parabolic bundles is non--empty for large topological invariants. The
proof is a generalization of the classical Cayley-Bacharach
construction to the setting of orbifold bundles. Our proof of {\it
  non-emptiness} and existence of components with smooth points gives
the same results for the Maruyama-Yokogawa space as well in the case
when $X$ is a surface. To the best of our knowledge the non-emptiness
of these moduli spaces have not been shown hitherto. In the sequel
(\cite{balasuman}) we also show the asymptotic irreducibility and
asymptotic normality of these spaces.

The moduli spaces are defined when some natural topological invariants
of the underlying objects are kept fixed. We also relate the
topological invariants that occur in (\cite{kron}, \cite{kron1}) with
natural invariants for parabolic bundles namely parabolic Chern
classes as defined in \cite{Biswaschern}. One observes that the
concept of an {\it action} (as defined in \cite{kron}) of a parabolic
bundle is precisely the {\it second parabolic Chern class}.  Moreover,
when we examine the Donaldson-Uhlenbeck compactification for these
moduli spaces, as observed by Kronheimer and Mrowka, the falling of
the instanton numbers is not perceived very precisely but what is seen
to drop in the boundary is the {\it second parabolic Chern class} or
equivalently the {\it action}.  Indeed, this is an exactly the
phenomenon in the usual Donaldson-Uhlenbeck compactification of stable
$SU(2)$-bundles on surfaces. For applications involving Donaldson
invariants arising from moduli of parabolic bundles should yield
topological invariants for the pair $(D,X)$ together with the
imbedding $D \hra X$ we refer the reader to \cite{kron}.

\medskip {\em Acknowledgment.}\, We are extremely grateful to
D.S.Nagaraj for his assistance and his invaluable comments and
suggestions. We thank C.S.Seshadri and S.Bandhopadyay for some useful
discussions.

\section{Preliminaries}
\subsection{\it {The category of bundles with parabolic structures}}

We rely heavily on the correspondence between the category of
parabolic bundles on $X$ and the category of $\G$--bundles on a
suitable Kawamata cover. This strategy has been employed in many
papers (for example \cite{Bi1}) but since we need its intricate
properties, most of which are scattered in a few papers of Biswas and
Seshadri, we recall them briefly. We stress only on those points which
are relevant to our purpose.

Let $D$ be an effective divisor on $X$. For a coherent sheaf
$E$ on $X$ the image of $E\bigotimes_{{\co}_X} {\co}_X(-D)$
in $E$ will be denoted by $E(-D)$. The following definition
of parabolic sheaf was introduced in \cite{MY}.

\medskip
{\bf Definition 2.3.}\, Let $E$ be a torsion-free
${\co}_X$--coherent sheaf
on $X$. A {\it quasi--parabolic} structure on
$E$ over $D$ is a filtration by ${\co}_X$--coherent
subsheaves
$$
E\, =\, F_1(E)\, \supset\, F_2(E)\, \supset\, \cdots
\,\supset\, F_l(E)\,\supset\, F_{l+1}(E)\,=\, E(-D)
$$
The integer $l$ is called the {\it length of the filtration}.
A {\it parabolic structure} is a quasi--parabolic structure,
as above, together with a system of {\it weights}
$\{{\a}_1,\cdots ,{\a}_l\}$ such that
$$
0\, \leq\,
{\a}_1\, < \, {\a}_2 < \, \cdots \, < \, {\a}_{l-1} \, < \,
{\a}_l \, < \, 1
$$
where the weight ${\a}_i$ corresponds to the subsheaf $F_i(E)$.
\medskip

We shall denote the parabolic sheaf defined above
by $(E,F_*,{\a}_*)$.
When there is no scope of confusion it will be denoted by $E_*$.

For a parabolic sheaf $(E,F_*, {\a}_*)$ define
the following filtration $\{E_t\}_{t\in \bR}$ of coherent
sheaves on $X$ parameterized by $\bR$:
\beqa
E_t \hspace{.1in} := \hspace{.1in} F_i(E)(-[t]D)
\eeqa
where $[t]$ is the integral part of $t$
and ${\a}_{i-1} < t - [t]
\leq {\a}_i$, with the convention that ${\a}_0 = {\a}_l -1 $
and ${\a}_{l+1} = 1$.

A {\it homomorphism} from the parabolic sheaf
$(E, F_*, {\a}_*)$ to another parabolic sheaf
$(E', F'_*, {\a}'_*)$ is a homomorphism from $E$ to $E'$
which sends any subsheaf $E_t$ into $E'_t$, where $t \in [0,1]$
and the filtration are as above.

If the underlying sheaf $E$ is locally free then $E_*$ will be called
a parabolic vector bundle. {\it In this section, all parabolic sheaves
  will be assumed to be parabolic vector bundles.}

\brem The notion of {\it parabolic degree} of a parabolic bundle $E_*$
of rank $r$ is defined as:
\beqa\label{parabdeg}
par_{deg}(E_*) := \int^{1}_{0} deg(E_t) dt + r.deg(D)
\eeqa
Similarly one may define $par_{\mu}(E_*) := par_{deg}(E_*)/r$. There is
a natural notion of parabolic subsheaf and given any subsheaf of $E$
there is a canonical parabolic structure that can be given to this
subsheaf. (cf \cite{MY} \cite{Bi1} for details)
\erem

\bdefe A parabolic sheaf $E_*$ is called {\em parabolic
  semistable (resp parabolic stable)} if for every parabolic subsheaf
$V_*$ of $E_*$ with $0 < rank(V_*) < rank(E_*)$, the following holds:
\beqa\label{parabss}
par_{\mu}(V_*) \leq par_{\mu}(E_*) ~~(resp. par_{\mu}(V_*) < par_{\mu}(E_*))
\eeqa
\edefe

\subsubsection{\it Some assumptions}\label{assumptions}
The class of parabolic vector bundles that are dealt with in the
present work satisfy certain constraints which will be explained now.
In a remark below, (see Remark \ref{weights}), we observe that these
constraints are not stringent in so far as the problem of moduli spaces
is concerned. 

\begin{enumerate}

\item  The first condition is that all parabolic divisors are
assumed to be {\it divisors with normal crossings}. In other words,
any parabolic divisor is assumed to be reduced, its each irreducible
component is smooth, and furthermore the irreducible components
intersect transversally.  

\item The second condition is that all the parabolic weights are {\it
    rational numbers}. 
  
\item The third and final condition states that on each component of
  the parabolic divisor the filtration is given by {\it subbundles}.
  The precise formulation of the last condition is given in
  (\cite{Bi1}, Assumptions 3.2 (1)\label{assumptions}). {\it
    Henceforth, all parabolic vector bundles will be assumed to
    satisfy the above three conditions.}

\end{enumerate}

\brem\label{weights} We remark that for the purpose of construction of
the moduli space of parabolic bundles the choice of rational weights
is not a serious constraint and we refer the reader to \cite[Remark
2.10]{MS} for more comments on this.  \erem

\bdefe A quasi--parabolic filtration on a sheaf $E$ can also be
defined by giving filtration by subsheaves of the restriction
$E|_{D}$ of the sheaf $E$ to each component of the parabolic divisor:
\[
E|_{D} = {\mathcal F}^1_{D}(E) \supset {\mathcal F}^2_{D}(E) \supset
\ldots \supset {\mathcal F}^l_{D}(E) \supset {\mathcal F}^{l+1}_{D}(E) = 0
\] 
together with a system of weights 
\[
0\, \leq\,
{\a}_1\, < \, {\a}_2 < \, \cdots \, < \, {\a}_{l-1} \, < \,
{\a}_l \, < \, 1
\]

\edefe

Let ${\rm PVect}(X,D)$ denote the category whose
objects are parabolic vector
bundles over $X$ with parabolic structure over the divisor
$D$ satisfying the above three conditions, and the
morphisms of the category are homomorphisms of parabolic
vector bundles (which was defined earlier).

The direct sum of two vector bundles with parabolic structures has an
obvious parabolic structure.  Evidently ${\rm PVect}(X,D)$ is closed
under the operation of taking direct sum. We remark that the category
${\rm PVect}(X,D)$ is an additive tensor category with the direct sum
and the parabolic tensor product operation. It is straight--forward to
check that ${\rm PVect}(X,D)$ is also closed under the operation of
taking the parabolic dual defined in \cite{Bi3Yo}.

For an integer $N\geq 2$, let ${\rm PVect}(X,D,N) \,
\subseteq \, {\rm PVect}(X,D)$
denote the subcategory consisting of all parabolic
vector bundles all of whose parabolic weights are multiples of
$1/N$. It is straight--forward to check that ${\rm PVect}(X,D,N)$
is closed under all the above
operations, namely parabolic tensor product, direct sum and
taking the parabolic dual.

\subsection{\it {The Kawamata Covering lemma}}

The ``Covering Lemma'' of Y. Kawamata
(Theorem 1.1.1 of \cite{KMM}, Theorem 17 of \cite{K}) says
that there is a connected smooth projective
variety $Y$ over $\bC$ and a Galois covering morphism
\beqa
p \,: \hspace{.1in} Y\hspace{.1in} \longrightarrow
\hspace{.1in} X 
\eeqa
such that the reduced divisor $\tD:= \,({p}^*D)_{red}$
is a normal crossing divisor on $Y$ and furthermore,
${p}^*D_i= k_iN.({p}^*D_{i})_{red}$,
where $k_i$, $1\leq i \leq c$,
are positive integers. Let $\G$ denote the Galois group
for the covering map $p$.

\subsection{\it{The category of $\G$--bundles}}\label{jumbo}
Let $\G \, \subseteq \, \mbox{Aut}(Y)$ be a finite subgroup
of the group of automorphisms of a connected smooth
projective variety $Y/{\bC}$. The natural
action of $\G$ on $Y$ is encoded in a morphism
$$
{\mu} \,:\hspace{.1in} \G \times Y \hspace{.1in}
\longrightarrow \hspace{.1in} Y
$$
Denote the projection of $\G\times Y$ to $Y$ by $p_2$. The
projection of $\G\times\G\times Y$ to the $i$--th factor will be
denoted by $q_i$. A {\it $\G$--linearized vector bundle} on $Y$ is a
vector bundle $V$ over $Y$ together with an isomorphism
$$
{\lambda} \,:\hspace{.1in} p^*_2 V \hspace{.1in}
\longrightarrow \hspace{.1in} {\mu}^*V
$$
over $\G\times Y$
such that the following diagram of vector bundles
over $\G\times\G\times Y$ is commutative:
\begin{center}
\begin{picture}(400,325)(0,20)
\put(25,290){\ctext{$q^*_3V$}}
\put(75,290){\vector(1,0){200}}
\put(425,290){\ctext{$(\mu \circ (q_2,q_3))^*V$}}
\put(325,250){\vector(0,-1){200}}
\put(75,265){\vector(1,-1){200}}
\put(325,25){\ctext{$(\mu\circ (m,Id_Y))^*V$}}
\put(50,150){\ctext{$(m\times Id_Y)^*\lambda$}}
\put(475,150){\ctext{$(Id_{\G}\times\mu)^*\lambda$}}
\put(175,320){\ctext{$(q_2,q_3)^*\lambda$}}
\end{picture}
\end{center}
where $m$ is the multiplication operation on $\G$.

The above definition of $\G$--linearization is equivalent
to giving isomorphisms of vector bundles
$$
{\bar g} \, : \, V \hspace{.1in} \longrightarrow
\hspace{.1in} (g^{-1})^*V
$$
for all $g \in \G$, satisfying the
condition that $\overline{gh}\,
= \, {\bar g}\circ {\bar h}$ for any $g,h \in \G$.

A {\it $\G$--homomorphism} between two $\G$--linearized vector bundles
is a homomorphism between the two underlying vector bundles which
commutes with the $\G$--linearizations. Clearly the tensor product of
two $\G$--linearized vector bundles admits a natural
$\G$--linearization; so does the dual of a $\G$--linearized vector
bundle. Let ${\rm Vect}_{\G}(Y)$ denote the additive tensor category
of $\G$--linearized vector bundles on $Y$ with morphisms being
$\G$--homomorphisms.

As before, ${\rm Vect}_{\G}(Y)$ denotes the category of
all $\G$--linearized vector bundles on $Y$.
The isotropy group of any point $y \in Y$, for the
action of $\G$ on $Y$, will be denoted by ${\G}_y$.

\subsection{On local types of $\G$--bundles}
Recall that since the $\G$--action on $Y$ is {\it properly
  discontinuous}, for each $y \in Y$, if $\G_y$ is the isotropy
subgroup at $y$, then there exists an {\it analytic} neighbourhood
$U_y \subset Y$ of $y$ which is $\G_y$-invariant and such that for
each $g \in G$, $g \cdot U_y \cap U_y \neq \emptyset$. 

\bdefe\label{bundlesfromreps} Let $\rho$ be a representation of $\G$
in $GL(r, \bC)$. Then $\G$--acts on the trivial bundle $Y \times
{\bC}^r$ by $(y,v) \lr (\gamma y, \rho(\gamma) v), y \in Y, v \in
{\bC}^r, \gamma \in \G$. Following \cite{S1} we call this
$\G$--bundle, the $\G$--bundle associated to the representation
$\rho$. \edefe

We then have the following equivariant local trivialisation lemma.

\blem\label{localtrivial} Let $E$ be a $\G$--bundle on $Y$ of rank
$r$. Let $y\in Y$ and let $\G_y$ be the isotropy subgroup of $\G$ at
$y$. Then there exists a $\G_y$-invariant analytic neighbourhood
$U_y$ of $y$ such that the $\G_y$--bundle $E|_{U_y}$ is associated to a
representation $\G_y \to GL(r)$ (in the sense of Def
\ref{bundlesfromreps}). \elem

\brem The above Lemma for $\G$--bundles with structure group $GL(r)$
can be found in \cite[Remark 2, page 162]{S1} and \cite{groth}.  Here
the key property that is used is that $U_y$ and $U_y/{\G_y}$ are Stein
spaces. This result, for the more general setting of arbitrary compact
groups $K$ instead of $\G$ and for general structure groups can be
found in \cite[Section 11]{heinzner}. \erem

\subsubsection{\it $\G$--bundles of fixed local type}\label{localtype}
We make some general observations on the local structure of
$\G$--bundles on the Kawamata cover defined in \eqref{jumbo}.

Let ${\rm Vect}^D_{\G}(Y,N)$ denote the subcategory of ${\rm
  Vect}_{\G}(Y)$ consisting of all $\G$--linearized vector bundles $W$
over $Y$ satisfying the following two conditions:

\begin{enumerate}
  
\item{} for a general point $y$ of an irreducible component of
  $(p^*D_i)_{red}$, the isotropy subgroup ${\G}_y$ is {\it
    cyclic} of order $|{\G}_y| = n_y$ which is a divisor of $N$; the
  action of the isotropy group ${\G}_y$ on the fiber $W_y$ is of order
  $N$, which is equivalent to the condition that for any $g \in
  {\G}_y$, the action of $g^{N}$ on $W_y$ is the trivial action;
  
\item{} In fact, the action is given by a representation $\rho_y$ of
  ${\G}_y$ given as follows:
  
  \beqa \rho_y(\zeta)= \left[
\begin {array}{cccc}
{z}^{{\alpha}_1}.I_1&&&0\\
&.&&\\
&&.&\\

0&&&{z}^{{\alpha}_l}.I_l
\end {array}
\right]\eeqa 

where 
\begin{itemize}
\item $\zeta$ is a generator of the group ${\G}_y$ and whose order
$n_y$ divides $N$ 

\item $\alpha_i = \frac{m_j}{N}$ and 
  
\item $I_j$ is the identity matrix of order $r_j$, where $r_j$ is the
  multiplicity of the weight $\alpha_j$. 
  
\item $z$ is an $n_y$-th root of unity. 
  
\item We have the relation $0 \leq m_1 < m_2 < ...< m_l \leq N - 1 $.

\end{itemize}
  
\item For a general point $y$ of an irreducible component of a
  ramification divisor for $p$ not contained in $(p^*D)_{red}$,
  the action of ${\G}_y$ on $W_y$ is the trivial action.
  
\item For a {\it special point} $y$ contained in $(p^*D)_{red}$, the
  isotropy subgroup $\G_y$ contains the cyclic group $\G_n$ of order
  $n$ determined by the irreducible component containing $y$. By the
  rigidity of representations of finite groups, the $\G_y$--module
  structure on $W_y$ (given by Lemma \ref{localtrivial}) when
  restricted to $\G_n \subset \G_y$ is of type $\tau$.
  
\item At special points $y$ of the ramification divisor for $p$ not
  contained in $(p^*D)_{red}$, the restriction of the
  representation to the generic isotropy is trivial.

\end{enumerate}

\bdefe Following Seshadri \cite[page 161]{S1} we call the $\G$--bundles
$E$ in ${\rm Vect}^D_{\G}(Y,N)$ {\sf bundles of fixed local orbifold
  type $\tau$}. \edefe

\brem The reason for calling it local type $\tau$ is that, for a
$\G$--bundle and a point $y$ the generic point of a divisor as above,
the structure of the representation defines the bundle $E_U$ for a
$\G_y$-invariant analytic neighbourhood in $Y$. Seshadri denoted the
collection of representations of the cyclic groups which define the
local isomorphism type over an analytic neighbourhood by the letter
$\tau$; note that the $\G$--bundle defines what is known as an
orbifold bundle. \erem

\brem We remark that this definition of $\G$--bundles of fixed local
type easily extends to $\G$--torsion--free sheaves since the local
action is specified only at the generic points of the ramification
divisor. \erem

We note that ${\rm Vect}^D_{\G}(Y,N)$ is also an additive tensor
category.

\subsubsection{\it {Parabolic bundles and $\G$--bundles}}\label{parabgamma}

In \cite{Bi1} an identification between the objects of ${\rm
  PVect}(X,D,N)$ and the objects of ${\rm Vect}^D_{\G}(Y,N)$ has been
constructed. Given a $\G$--homomorphism between two $\G$--linearized
vector bundles, there is a naturally associated homomorphisms between
the corresponding vector bundles, and this identifies, in a bijective
fashion, the space of all $\G$--homomorphisms between two objects of
${\rm Vect}^D_{\G}(Y,N)$ and the space of all homomorphisms between
the corresponding objects of ${\rm PVect}(X,D,N)$.  An equivalence
between the two additive tensor categories, namely ${\rm
  PVect}(X,D,N)$ and ${\rm Vect}^D_{\G}(Y,N)$, is obtained this way.
Since the description of this identification is already given in
\cite{Bi1}, and \cite{Bi3}, it will not be repeated here.



We observe that an earlier assertion that the parabolic tensor product
operation enjoys all the abstract properties of the usual tensor
product operation of vector bundles, is a consequence of the fact that
the above equivalence of categories indeed preserves the tensor
product operation.

The above equivalence of categories has the further property that it
takes the parabolic dual of a parabolic vector bundle to the usual
dual of the corresponding $\G$--linearized vector bundle.

Let $W\, \in \,{\rm Vect}^D_{\G}(Y,N)$ be the $\G$--linearized vector
bundle of rank $n$ on $Y$ that corresponds to the given parabolic
vector bundle $E_*$.  The fiber bundle
$$
\pi \, :\hspace{.1in} P \hspace{.1in} \longrightarrow
\hspace{.1in} Y
$$
whose fiber ${\pi}^{-1}(y)$ is the space of all $\bC$--linear
isomorphisms from ${\bC}^n$ to the fiber $W_y$, has a
the structure of a $(\G ,GL(n, {\bC}))$--bundle over $Y$.

\bdefe A $\G$-linearized vector bundle $E$ over $Y$ is called
$\G$-{\it semistable} (resp. $\G$-{\it stable}) if for any proper
nonzero coherent subsheaf $F \subset E$, invariant under the action of
$\G$ and with $E/F$ being torsionfree, the following inequality is
valid: 
\beqa\label{gammass} {\mu}(F) \leq {\mu}(E)~~ (resp. {\mu}(F) <
{\mu}(E)) 
\eeqa 
where the slope is as usual ${\mu}(E) = deg(E)/r$ and $deg(E)$ is
computed with respect to the $\G$--linearised very ample divisor
$\Theta$ on $Y$.

The $\G$-linearized vector bundle $E$ is called $\G$-{\it polystable}
if it is a direct sum of $\G$-stable vector bundles of same slope.

\edefe

\brem\label{piparab} The above correspondence between {\it parabolic
  bundles} on $X$ and $\G$--bundles on $Y$ preserves the semistable
(resp. stable) objects as well, where {\em parabolic semistability} is
as in \eqref{parabss}. (cf \cite{Bi1}) \erem

\brem \label{gammasemistability} We remark that it is not hard to
check that for $\G$--bundles, $\G$--semistability (resp.
$\G$--polystability) is the same as usual semistability (resp.
polystability).  This can be seen from the fact that the {\em top term
  of the Harder-Narasimhan filtration} (resp. {\em the socle}) are {\it
  canonical} and hence invariant under the action of $\G$. But we note
that a $\G$--stable bundle need not be $\G$--stable, as can be seen by
taking a direct sum of $\G$-translates of a line bundle. \erem

\brem\label{gammacohomology} 

We make some key observations in this remark where we also note the
essential nature of assumptions of characteristic zero base fields.
\begin{enumerate}
\item The notion of $\G$--cohomology for $\G$--sheaves on $Y$ has been
  constructed and dealt with in great detail in \cite{grothendieck}.
  These can be realised as higher derived functors of the {\it
    $\G$--fixed points}--sub-functor $({{H}^{0}})^{\G}$ of
  the section functor ${H}^{0}$. (We use this notation to
  avoid $\G^\G$, because we have denoted the finite group by the
  letter $\G$!).
  
  We note immediately that since we work over fields of characteristic
  zero, the sub-functor $({{H}^{0}})^{\G} \subset {{H}^{0}}$ is in
  fact a direct summand (by averaging operation).  Hence, we see
  immediately that the higher derived functors of the
  functor$({{H}^{0}})^{\G}$ are all sub objects of the
  derived functors of ${{H}^{0}}$.

\item When we work with a Kawamata cover as in our case, then we have
  the following relation between the $\G$--cohomology and the usual
  cohomology on $Y/{\G} = X$:

\[
H^{i}_{\G}(Y, \mathcal F) = H^i(X, p^{\G}_{*}(\mathcal F))
\]
$\forall i$.

\end{enumerate}

\erem

\subsubsection{\it $\G$--bundles and orbifold bundles}

We make a few general remarks on the advantages of working with a
Kawamata cover $Y$ and $\G$--bundles on $Y$ over working with {\it
  orbifold bundles} or $V$--bundles over $V$--manifolds. Locally,
these two notions can be completely identified but for any global
construction such as the one which we intend doing, namely a moduli
construction, working with a Kawamata cover albeit {\it non-canonical},
has obvious advantages since it immediately allows us to work with a
certain ``Quot'' scheme over $Y$. To recover the moduli of parabolic
bundles with fixed quasi parabolic structure, we then simply use the
{\it functorial equivalence of parabolic bundles and $\G$--bundles of
  fixed local type}. 

\subsubsection{\it $\G$--line bundles and parabolic line bundles}\label{plb} 
A $\G$ line bundle on $Y$ is a line bundle $L$ on $Y$ together with a
lift of action $\G$ . The $\G$ line bundle gives a $\G$ invariant line
bundle $L^{\G}$ on $X$.  Let $D$ be a divisor of normal crossing on
$X$.  Let $D=\sum_{i=1}^{d} D_i$ be a decomposition into irreducible
components.  A parabolic line bundle on $(X, D)$ is a pair of the form
$(M,{\beta_1,...,\beta_i,...,\beta_d})$ where $M$ is a holomorphic
line bundle on $X$ and $0\le \beta_i < 1$ is a real number. When we
start from a $\G$ line bundle on $Y$ we get a pair
$(L^{\G},{\beta_1,...,\beta_i,...,\beta_d})$ where $\beta_i$ is a
rational number and it can be written as $\beta_i=m_i/N.$ Let
$\tilde{D_i} = (p^*D_i)_{red}$. Then by following \cite[Section
2b]{Bi4} we have $L=p^*(L^\G)\otimes \m O_Y(\sum_{i=1}^d k_i
m_i\tilde{D}_i)$

\brem In our situation, by choice we work with a single weight when we
consider $\G$--line bundles of fixed local type $\tau$ although this
may not be absolutely essential.\erem

\subsubsection{\it Serre duality for $\G$--line bundles of fixed local
  type}\label{serredual} 

\bdefe\label{taulinebundles} By a line bundle $L$ of fixed local type
$\tau$ we mean a {\it parabolic line bundle} $(L, \alpha_1, \alpha_2,
\ldots, \alpha_d)$, where $\alpha_i = \alpha \forall i$. In other
words, locally, the generic isotropy on the irreducible components of
the inverse image of the parabolic divisor acts by a single character
namely $\alpha$. We will write $L^{(\alpha)}$ to specify the character.
\edefe

Let $L = L^{(\alpha)}$ be a $\G$ line bundle on $Y$ of type $\tau$.
Then by \ref{plb}, one knows that $ L = p^*({p^{\G}_*}(L))\otimes \m
O_Y(\sum k_i m_i \tilde{D}_i) $ where all the $m_i$ can be assumed to
be equal to $m$ since we have a single weight $\alpha$. Then if $M =
M^{(\alpha)}$ is another $\G$--line bundle with the same local
character type we have:

\beqa\label{gammaiso}
({p^{\G}_*}(L^*\otimes M)) = ({p^{\G}_*}(L)^*\otimes ({p^{\G}_*}(M))
\eeqa 

Consider the canonical bundles $K_X$ of $X$ and define the
$\G$--bundle $K_Y^{(\alpha)}$ as follows:
\beqa\label{twistedcanonical}
K_Y^{(\alpha)} = p^*(K_X) \otimes \m O_Y(\sum k_i {\tilde{D}}_i) m)
\eeqa
Then, we see as above that ${p^{\G}_*}(K_Y^{(\alpha)}) = K_X$. We then
have the following duality for $\G$--line bundles of type $\tau$:

\blem\label{serreduality} For $\G$--line bundles $L$ of type $\tau$,
with local character $\alpha$, the $\G$--line bundle $K_Y^{(\alpha)}$
is the {\em dualising sheaf}. In other words, we have a canonical
isomorphism:
\[
H^{i}_{\G}(Y, L^{*} \otimes K_Y^{(\alpha)}) \simeq H^{n-i}_{\G}(Y,
L)^{*}
\]
for all $i$. We have made this statement for $\G$--varieties $Y$ of
any dimension.
\elem

\pf The proof is straightforward, but we give it for the sake of
completeness. Recall the relationship between the $\G$--cohomology on
$Y$ and the usual cohomology on $X$ (Remark \ref{gammacohomology}). We
have the following isomorphism (using \ref{gammaiso}):
\[
H^{i}_{\G}(Y, L^{*} \otimes K_Y^{(\alpha)}) \simeq H^{i}(X,
{p^{\G}_*}(L^*\otimes K_Y^{(\alpha)}) \simeq H^{i}(X, {p^{\G}_*}(L)^*\otimes
({p^{\G}_*}(K_Y^{(\alpha)})) 
\]
Using ${p^{\G}_*}(K_Y^{(\alpha)}) = K_X$ we then conclude from the following
isomorphism:
\[
\simeq H^{i}(X, {p^{\G}_*}(L)^* \otimes K_X)\simeq H^{n-i}(X,
{p^{\G}_*}(L))^{*} \simeq H^{n-i}_{\G}(Y, L)^{*}
\]
where we use the usual Serre duality on $X$.
\enpf

\section{Towards the construction}

\subsubsection{On determinant line bundles}\label{detlinebundles}

We briefly recall the basic definitions for the convenience of the
reader. Let $Y$ be an irreducible smooth projective variety equipped
with a very ample ${\mathcal O}_Y(1)$. Let $K(Y)$ be the Grothendieck
algebra of classes of coherent sheaves. Let $\theta$ be the class in
$K(Y)$ of the structure sheaf ${\mathcal O}_{\Theta}$ of a hyperplane
section $\Theta \subset Y$. This algebra is equipped with a quadratic
form $q : u \mapsto {\chi}(u^{2})$. This form is calculated in terms
of the rank and the Chern classes of $u$. For example, if $Y$ is a
smooth projective surface, and if $u \in K(Y)$ is of rank $r$, and the
Euler characteristic $\chi$, we have
\[
q(u) = 2 r {\chi} + c_1^{2} - r^{2} {\chi}({\mathcal O}_Y)
\]
The kernel $ker (q)$ comprises of the classes which are {\it numerically
  equivalent to zero}.  We work with the quotient:
\[
K_{num}(Y) = K(Y)/{ker (q)}
\]
For a smooth projective surface $Y$, $K_{num}(Y) \simeq {\mathbb Z} \times
H^{2}(Y, \mathbb Z) \times {\mathbb Z}$. and this isomorphism is by giving
$(r,c_1,\chi)$.  

Recall that if $\cf$ is a flat family of coherent sheaves on $Y$
parametrised by a scheme $S$, then $\cf$ defines an element $[\cf] \in
K^{0}( S \times Y)$, the Grothendieck group of $S \times Y$ generated
by locally free sheaves.  We may then define the homomorphism from the
Grothendieck group of coherent sheaves on $Y$ given by:
\[
{\lambda}_{\cf} : K(Y) \lr Pic(S).
\]
as follows: For $u \in K(Y)$, ${\lambda}_{\cf}(u)= det (pr_{1!} ({\cf}
\cdot pr_{2}^{*}(u))$, where ${\cf} \cdot pr_{2}^{*}(u)$ is the product
in $K(S \times Y)$ and $pr_{1!} : K^0(S \times Y) \rightarrow K^0(S)$
associates to each class $u$ the class $\sum_i(-1)^i
R^ipr_{1*}({u})$.

We observe that this has a collection of functorial properties for
which we refer to (\cite{H} page 179).

Let $Y$ be a smooth projective surface. Fix a class $c \in
K_{num}(Y)$, i.e the rank $r$, the first Chern class $c_1 = {\co}_Y$
and the Euler characteristic $\chi$. This in particular fixes $c_2$ as
well.  Fix also the very ample divisor $\Theta$ on $Y$ and a base
point $x \in Y$. Let $\theta = [{\mathcal O}_{\Theta}] \in K(Y)$.
Define for each $i$:

\beqa\label{ui} u_i(c) := - r \cdot {\theta}^i + {\chi}(c \cdot
{\theta}^i) \cdot [{\mathcal O}_x] \eeqa 

(cf \cite[page 183]{H}).

\subsection{\it{Projective $\G$--frame bundle}}\label{frame}
We make some general remarks on the general construction of
$\G$--frame bundle associated to a $\G$--vector bundle. This is a
generalization of the classical frame bundle construction but will be
needed in the construction of the moduli space. Let $Y$ be a scheme of
finite type with a {\it trivial} $\G$--action.  Let $F$ be a
$\G$--locally free $\m O_Y$ module of rank $r$ and assume that each
fibre $F_y$ is a $\G$--module and the $\G$--module structures are
isomorphic at different points. Let $W$ be a finite dimensional vector
space of dimension $r$ which is a $\G$--module isomorphic to the
$\G$--module $F_y$ for any $y \in Y$.  Denote by $\m O_Y(W)$ the
trivial rank $r$ sheaf modelled by $W$. With this added structure, we
have a canonical group namely, $H = Aut_{\G}(W) \subset GL(W)$, which
acts on $\m O_Y(W)$ by automorphisms which preserve the
$\G$--structure.

Let $\mathbb Hom_\G(\m O_Y(W), F):= Spec(S^*(\mathcal Hom_\G(\m
O_Y(W), F)))^* \rightarrow Y$ be the geometric $\G$--vector bundle that
parameterises all $\G$--homomorphisms from $\m O_Y(W)$ to $F$.  Let
${\Phi(F)} := \mathbb Isom_\G(\m O_Y(W), F) \subset \mathbb Hom_\G(\m
O_Y(W), F)$ be the open subscheme which parameterises all
$\G$--isomorphisms and let $\pi:{\Phi(F)} \rightarrow Y$ denote the
canonical projection. 

Then we observe that $H$ acts on $\Phi(F)$ by composition and $\pi$ is
a principal bundle with structure group $H$. Indeed, the
$\G$--structure on $F$ gives a natural reduction of structure group of
the frame bundle associated to $F$ (which by the usual construction is
a principal $GL(W)$--bundle).

Similarly, if $PH$ is the image of $H \subset GL(W)$ in $PGL(W)$, then
one can construct projective $PH$--bundle by taking image of $\Phi(F)$
in $Proj(S^*(\mathcal Hom_\G(\m O_Y(W), F)^*)) $. We term the image of
$\Phi(F)$ the {\it projective $\G$--frame bundle} over $Y$ associated
to the $\G$--bundle $F$.


\subsection{\it{The determinant line bundle}}

The aim of this section is to construct a line bundle on the Quot
scheme which parametrises the objects we need. This will be a natural
{\it determinantal bundle} as in the Donaldson construction. 

Recall that our aim is to construct the moduli space of
$\mu$--semistable bundles with $\G$--structure and the notion of
$\mu$--semistability in the higher dimensional setting (in our case
the surface $Y$) is not a GIT notion; in fact, the GIT semistable will
be the Gieseker semistable bundles.

Since $\G$--semi stability is the same as usual semistability for
torsion free sheaves (cf Remark \ref{gammasemistability}) we observe
that the family of $\G$--semistable sheaves with fixed Hilbert
polynomial is bounded (Thm. 3.3.7 \cite{H}).

Let $\E$ be a torsion free $\G$--coherent sheaf over a smooth
projective surface $Y$, of rank $r$ and $P$ be any polynomial in
$\mathbb Q[z]$. $\text{Quot}(\m E, P)$ be the Quot scheme which
parametrises all quotients of $\m E$ with fixed Hilbert polynomial
$P$. Let $\m F$ denote the universal quotient sheaf of $\m
O_{Quot(\E,P)} \otimes \E$ on $Y \times Quot(\E,P)$. {\sf Let $Q$
  denote the subscheme of $Quot(\E,P)$} whose closed points correspond
to torsion--free sheaves with fixed topological data $(c_1,c_2,r)$
(note that fixing Hilbert polynomial for a family of sheaves gives
only finitely many choices for the triplets $(c_1,c_2,r)$) and $\m F
\mid_{Q \times Y}$ be universal quotient sheaf on $ Q \times Y$. Let
$L$ be the determinantal line bundle ${\lambda}_{\m F}(u)$. Since $\G$
is acting on $\E$ and $Y$, $\G$ acts on $Q$ in the natural manner:
\[
\xymatrix{
\E \ar[r]^{[q]} & \m F_q  \\ 
\E \ar[u]_{\gamma^*} \ar[ur]\\
}
\]

where $\gamma^*$ is the canonical pull back.  Let $Q^\G \subset Q$ be
the set of all $\G$--invariant points of $Q$ which is a nonempty
subset (!), and by following \cite{S1} it gets a {\it closed subscheme
  structure}.

Let $P_c(m) = {\chi}(c(m))$ be the Hilbert polynomial associated to
the fixed class $c \in K_{num}(Y)$, where $c(m) := c \cdot [{\mathcal
  O}_Y(m)]$.  Let $\E = V \otimes \mathcal O_{Y}(-m)$ where $V$ is a
vector space of dimension $P_c(m)$. We choose $m$ large enough so that
all quotients are $m$ regular(i.e. higher cohomology group $H^i(Y,\m
F_q(m-i))$ vanishes for all $i \ge 1$ and for all quotients $\m F_q$
of $\Ha$).

\begin{note}\label{fixeddet} 
  Let $P = P_c(m)$ and let $Q = \text{Quot}(\m E, P)$. Let $Q^{\G}$
  denote the closed subscheme of $\G$--fixed points.  Let $\m R
  \subset Q$ (resp $\Rss \subset Q^\G$) be the locally closed
  subscheme of all $\mu$--semistable quotients (resp
  $(\G,\mu)$--semistable quotients) of $\m E$ with fixed topological
  data $(r,c_1,c_2)$ {\it and fixed determinant $\mathcal Q$}. We
  observe that giving the topological data is giving a class $c \in
  K_{num}(Y)$. \end{note}

Because of $m$--regularity we have $V \simeq H^0(\m F_q(m)) \simeq
\mathrm k^{P_c(m)}$. The group $Aut(V)$ acts naturally on the scheme
$Q$.

\begin{note} Let us denote by $G$ the group $SL(V)$ and by $H$ the
  subgroup $Aut_{\G}(V) \cap G$ i.e the subgroup of $G$ which are
  $\G$--automorphisms as well. We will use this notation through this
  entire paper. \end{note}

\brem\label{connectedreductive} The group $Aut_\G(V)$ is a direct
product of full linear groups and in particular connected and
reductive. The group $H$ is also therefore connected and reductive To
see this, observe that we can decompose $V$ as a $\G$--module into its
isotypical decomposition.  This decomposition gives the choice of a
torus in $SL(V)$ and the group $H$ is the centraliser of this torus;
indeed, $H$ is the Levi subgroup associated to the parabolic subgroup
given by the decomposition. This implies that $H$ is connected and
reductive. The group $Aut_\G(V)$ is similarly the Levi subgroup in the
bigger group $GL(V) = Aut(V)$ \erem

The group $H$ (resp $G$) acts on the scheme $\Rss$ (resp $\m R$) by
automorphisms. The universal quotient $\m F$ allows us to construct a
$G$--linearised line bundle $\m N$ on $\m R$ given as follows:
\[
{\m N} : = {\lambda_{\cf}(u_1(c))}
\]
where $u_i(c)$ is defined as in \eqref{ui}. Denote by $\m M$ the
restriction of this line bundle to $\m R^\G$.  That is:
\beqa\label{mM} \m M = \m N|_{\m R^\G} \eeqa

Let $\Rss(D,N)$ be the subset $\Rss$ consisting of $\G$--torsion-free
sheaves of fixed local type. 

\brem By the rigidity of representation of finite groups, it
follows that $\Rss(D,N)$ is both {\it open and closed} in $\Rss$.
Moreover, it is easily seen that $\Rss(D,N)$ is also invariant under
the action of $H$.\erem

\brem By definition, the line bundle $\m M$ comes with a canonical
$H$--linearisation.
\erem

Then we have the following:

\blem\label{huy} (\cite [Lemma 8.2.4]{H})
\begin{enumerate}
\item[1.]If $s \in \Rss$ is a point such that for a general high
  degree $\G$--invariant curve $C$, $ \m F_s \mid_C$ is semistable
  then there exists an integer $N > 0$ and an $H$--invariant
  section $\tilde \sigma \in H^0(\Rss,\m M^N)^{H}$ such that
  $\tilde \sigma(s) \ne 0$.
  
\item[2.] If $s_1$ and $s_2$ are two points in $\Rss$ such that for a
  general high degree $\G$--invariant curve $C$, $\m F_{s_1} \mid_C$
  and $\m F_{s_2} \mid_C$ are both semistable but not $S$--equivalent
  or one of them is semistable but other is not then there is a
  $H$--invariant section $\tilde \sigma$, in some tensor power of $\m
  M$ which separates these two points (i.e $\tilde \sigma(s_1)=0$ but
  $\tilde \sigma(s_2) \ne 0$).
\end{enumerate}
\elem

{\it Proof:} The proof (following ideas from Le Potier \cite{Le}) is
largely following the exposition in Huybrechts-Lehn(\cite{H}), But we
give all the main steps in the argument even at the risk of
repetition. This is because there are certain distinctive points in
this setting which needs to be highlighted, especially those relating
to the projective $\G$--frame bundle and the morphism to the quot
scheme of $\G$--bundles on a curve. In a sense these are precisely the
points which distinguish the possible $\G$--structures on a given
semistable bundle.
 
Since $\G$--semistability is same as usual semistability, one gets a
general high degree smooth curve $C \in \mid a{\Theta} \mid^\G, a \gg 0$,
such that, $\m F \mid_{\Rss \times C}$ produces a family of
generically semistable sheaves on $C$ with fixed topological data
$(r,\m Q \mid_C)$. Recall that $\m Q$ is the fixed determinant for
objects in $\Rss$ (see \eqref{fixeddet}). The fact that it is a
generic family of semistable sheaves on $C$ is because of openness of
semistability property (cf for example \cite{NR}). Let $U$ be a
nonempty open subset of $\Rss$ such that $\m F \mid_{U \times C}$
is a flat family of semistable sheaves on $C$.

Recall that we have fixed a class $c \in K^\G_{num}(Y)$. Let $c
\mid_C$ be its pull-back (or restriction) in $K^\G_{num}(C)$. Note
that $c\mid_C$ is completely determined by its rank $r$ and the line
bundle $\m Q|_C$.

Recall that $P_c(m) = {\chi}(c(m))$ is the Hilbert polynomial
associated to the fixed class $c \in K^\G_{num}(Y)$, where $c(m) := c
\cdot [{\mathcal O}_Y(m)]$. Let $P'(n) := P_{c\mid_C}(n)$. Then, by
computing the Euler characteristic from the exact sequence of sheaves
obtained by {\it restriction to the curve $C$}, we see that $P'$ is
given by the equation $P'(n) = P_c(n) - P_c(n-a)$, since $C \in \mid
a{\Theta} \mid^\G$.

Let $\m H' = \m O_C(-m')^{P'(m')}$ and $Q_C^\G \subset
\text{Quot}^\G_C(\m H',P')$ be the closed subset of quotients with
determinant $\m Q\mid_C$. Observe that $\m H'$ can be identified with
$W \otimes \m O_C(-m')$, where $W$ is a vector space of dimension
$P'(m')$.

Denote by $G_1$ the group $SL(W)$ and by $H_1$ the subgroup of $G_1$ given by:
\[
H_1 = G_1 \cap Aut_{\G}(W)
\]
As remarked earlier (Remark \ref{connectedreductive}), the group $H_1$
is also {\it connected and reductive}.

We also have a natural $H_1$--action on $Q_C^\G$ by automorphisms.

Let $\m O_{{Q_C}^\G} \otimes \m H' \twoheadrightarrow \tilde {\m F}'$
be the universal quotient and $L_C = \lambda_{\tilde {\m F}'}(u_0(c\mid_C)$
(see \eqref{ui} for the definition of $u_0(c)$). 

{\it One can check that $L_C \cong det(p_{Q_C^\G*}(\tilde {\m F}'))$
}. If $m'$ is sufficiently large the following holds:

\begin{enumerate}
\item{} Given a point $[q:\m H' \twoheadrightarrow {\tilde {\m F}'}_q] \in
  Q_C^\G$, the following assertions are equivalent:
\begin{enumerate}\label{frameb}
\item{} ${\tilde {\m F}'}_q$ is $\G$--semistable sheaf
  and $W \simeq H^0(C,{\tilde {\m F}'}_q(m'))$
\item{} $[q]$ is a semistable point in $Q_C^\G$ for the action of
  $H_1$ with respect to the linearization of $L_C$, i.e, there is an
  integer $\nu$ and a $H_1$--invariant section $ \sigma \in
  H^0(C,L_C^{\nu})^{H_1}$ such that $\sigma([q]) \ne 0$.
\end{enumerate}
\item{} Two points $[q_i:\m H' \rightarrow {\tilde {\m
      F}'}_{q_i}];i=1,2$ are separated by $H_1$--invariant sections if
  and only if either both are semistable points but ${\tilde {\m
      F}'}_{q_1}$ and ${\tilde {\m F}'}_{q_2}$ are not $S$--equivalent
  or else, one of them is semistable and other is not semistable.
\item{} $\tilde
{\m F} := \m F  \mid_{\Rss \times C}$ is $m'$ 
regular with respect to $\Rss$. 

\end{enumerate}
Note that $ p_*(\tilde {\m F } (m'))$ is a $\G$--locally free $\m
O_{\Rss}$ sheaf of rank $P'(m')$. The group $H_1$ acts on $Q_C^\G$.
Let $\pi:\tilde {\Rss} \rightarrow \Rss$ be the associated
$PH_1$--bundle, i.e the {\it projective $\G$--frame bundle} (by
\eqref{frameb} above, the conditions required in \eqref{frame} hold
good here). From the $H$--action on $\Rss$, we see that $\tilde
{\Rss}$ gets an $H$--action as well.

The projective $\G$--frame bundle $ \tilde {\Rss} $ parametrises a
quotient $\m O_{\tilde {\Rss}} \otimes \m H' \lr \pi^* \tilde {\m F}
\otimes \m O_{\pi}(1)$. So it gives rise to $H_1$--equivariant morphism
$\phi_{{\tilde {\m F}}}: \tilde {\Rss} \longrightarrow Q^\G_C$. We
note that $\tilde {\Rss}$ also carries an $H$--action on it induced
from $\Rss$.  So $\tilde {\Rss}$ carries an $(H_1 \times H)$--action.
So one gets the following diagram

\[
\xymatrix{
\tilde {R}^{\G } \ar[r]^{\phi_{\tilde{\m F} }}
\ar[d]^\pi & Q^\G_C \\
R^{\G } \\
}
\]

We now use the computations involving determinant bundles in
\cite[8.2]{H} and the functoriality of the determinant bundle and note
the fact that all the families involved which are defined over the
schemes $\Rss$ and $Q^\G_C$, are just the pull--backs of the ones on
the usual quot scheme. It therefore follows that the relation obtained
in \cite[8.2]{H} hold verbatim over the projective $\G$--frame bundle
$\tilde {\Rss}$ as well.

We note that, since the projective $\G$--frame bundle $\tilde {\Rss}$
is the reduction of structure group of the usual projective frame
bundle over $\m R$ restricted to $\Rss$, $\tilde {\Rss}$ is a closed
subscheme of $\tilde {\m R}$ over $\Rss$. Thus, if $\m M$ is as in
\eqref{mM}, we have

\[
\phi_{\tilde {\m F}}^{*}(L_C)^{deg(C)}~~\simeq~~\pi^{*}(\m M)^{a^2
  deg(Y)}
\] 

If $s$ is a $H_1$--invariant section of $L_C^{\nu deg(C)}$ for some
$\nu > 0$, then $\phi_{\tilde {\m F}}^*(s)$ is a $(H_1 \times
H)$--invariant section i.e an element of $H^0({\tilde
  {\Rss}},\phi_{\tilde {\m F}}^*(L_C)^{\nu deg(C)})^{H_1 \times H} =
H^0({\tilde {\Rss}},\pi^{*}(\m M)^{\nu a^2 deg(Y)})^{H_1 \times H}$.

Since $\pi: \tilde {\Rss} \rightarrow {\Rss}$ is a principal
$PH_1$--bundle, the section $\phi_{\tilde {\m F}}^*(s)$ will descend
to give an element in $H^0(\Rss,\m M^{\nu a^2 deg(Y)})^{H}$. In
other words, for each $\nu > 0$, we get a linear (injective) map:
\[
s_{\m F}:H^0(Q^\G_C, L_C^{\nu deg(C)})^{H_1} \rightarrow
H^0(\Rss,\m M^{\nu a^2 deg(Y)})^{H}
\]

Now let ${\m F}_q$ be a point in $\Rss$, i.e a $\G$--semistable
torsion free sheaf. By the Orbifold Mehta--Ramanathan restriction
theorem (Theorem \ref{orbimehtaramanathan}) it follows that there
exists a curve $C$ as above such that the restriction ${\m F}_q|{_C}$
is in $Q^\G_C$. Hence, by the usual GIT and Seshadri's theorem, there
exists a section $s \in H^0(Q^\G_C, L_C^{\nu deg(C)})^{H_1}$ for
some $\nu > 0$ which is {\it non--zero} at the point ${\m F}_q|{_C}$.

Following the map $s_{\m F}$ we get a section in $H^0(\Rss,\m M^{\nu
  a^2 deg(Y)})^{H}$ which is non--zero at $\m F_q$ proving the
lemma.

\begin{flushright} {\it q.e.d} \end{flushright}

We have the following immediate corollary from the first part of Lemma
\ref{huy}:

\bcor\label{lepot} There exists an integer $\nu > 0$ such that the
line bundle ${\m M}^{\nu}$ on $\Rss$ is generated by $H$--invariant
global sections.  \ecor

\section{Donaldson--Uhlenbeck compactification}

The aim of this section is to construct a {\it reduced algebraic
  scheme} i.e a {\it variety}, which is projective and whose points
give the analogue of the underlying set of points of the
Donaldson--Uhlenbeck compactification for $\G$--bundles on a smooth
projective algebraic surface with a $\G$--action. This, in conjunction
with the Kawamata covering lemma and the general ({\it parabolic
  bundles)--($\G$--bundles}) correspondence would enable us to
construct a projective variety whose underlying set of points
parametrise the natural analogue of Donaldson--Uhlenbeck
compactification of the moduli space of {\it $\mu$--stable parabolic
  bundles} on a surface $X$ with parabolic structure on a divisor with
normal crossings. We also describe the {\it boundary points} of the
compactification in terms of $\G$--bundles and $0$--cycles on the
surface $Y$ (and as a consequence on $X$ as well).

Since $\Rss$ is a quasi--projective scheme and since $\m M$ is $H$--{\it
  semi--ample}, there exists a finite dimensional vector space $A
\subset A_\nu:= H^0(\Rss,\m M^{\nu})^{H}$ that generates $\m M^\nu$; of
course, there is nothing canonical in the choice of $A$.

Let morphism $\phi_A: \Rss \rightarrow \mathbb P(A)$ be the induced
$H$--invariant morphism defined by the sections in $A$. 

But because of non-uniqueness of $A$ a different choices of subspace of
invariant sections gives rise to a different map $\phi_{A'}$ to a
different projective space $\mathbb P(A')$.

\bdefe We denote the by $M_A$ the {\sf schematic image} $\phi_A(\Rss)$
with the {\sf canonical reduced scheme structure}.  \edefe

\brem By the following result which may be titled $H$--properness, the
variety $M_A$ is {\it proper} and hence because of its
quasi--projectivity it is a projective variety. We note that we use the
term {\it variety} in a more general sense of an {\it reduced
  algebraic scheme of finite type} which need not be irreducible.  So
in what follows we will be working with the $\mathbb C$--valued points
of $M_A$.

\erem

\bprop If $T$ is a separated scheme of finite type over $k$, and if
$\phi: \m R^{\mu ss} \lr T$ is an $Sl(V)$ invariant morphism then
image of $\phi$ is proper over $k$.  \eprop 

\brem This is a consequence of the Langton type semistable reduction
theorem for $\G$--torsion free sheaves which we have shown in the
Appendix and some general schematic methods (cf \cite [Prop 8.2.5]{H}
for details). \erem

Let $A_\nu$ denote the vector space $H^0(\Rss,\m M^\nu)^H$, $\nu
\in \mathbb Z^{+};$ and Let $A \subset A_\nu$ be a finite dimensional
vector space which generates $\m M^\nu$. 

For any $d \ge 1$, let $A^d$ be the image of the canonical
multiplication map $f_d:A \otimes,\cdots,\otimes A (d--times)
\rightarrow A_{d\nu}$; in particular $A^1 = A$.

Let $A'$ be any finite dimensional vector subspace of $A_{dN}$
containing $A^d$. Then clearly the line bundle $\m M^{d\nu}$ is also
globally generated by $H$--invariant sections coming from the subspace
$A'$ and this is so for any $d \ge 0$.  

So we have $A \rightarrow A^d \subset A' $, and hence
a commutative diagram
\[
\xymatrix{
M_{A'} \ar[r]^{\pi_{A'/A}} & M_A \\
\Rss \ar[u]_{\phi_{A'}} \ar[ru]_{\phi_A}
}
\]
Since $M_A$ and $M_A'$ are both projective, the map $\pi_{A'/A}$ is a
finite map (pull-back of ample remains ample). So if we fix a $A$ as
above we get an inverse system (indexed by the $d \ge 1$) of
projective varieties $(M_{A'},\pi_{A'/A})$ dominated by the finite
type scheme $\Rss$.

\[
\xymatrix{
& \Rss \ar[dl] \ar[d]_{\phi_{A'}} \ar[dr]^{\phi_{A}} \\
\cdots \ar[r] & M_{A'} \ar[r]_{\pi_{A'/A}} & M_{A}  
}
\]

Hence the inverse limit of the system $(M_{A'},\phi_{*})$ is in fact one of
the $M_{A'}$'s where $A'$ is a finite dimensional subspace of
$H^0(\Rss,\m M^{n})^{H}$ which generates $\m M^n$.  

\bdefe We denote this {\it inverse limit variety} $M_{\G}$ and let
$\phi:\Rss \rightarrow M_{\G}$ be the canonical morphism induced by the
invariant sections coming from the subspace $A'$ associated to the
inverse limit.\edefe

\brem We will show that the moduli space of isomorphism classes
$(\G,\mu)$--stable locally free sheaves of fixed type $\tau$ and {\it
  fixed determinant $\m Q$} will be a subvariety of $M_\G$. This will
allow us to take the closure of the moduli space of stable bundle in
$M_\G$ and give it the reduced scheme structure.\erem

\brem The underlying set of points of this projective variety , namely
the closure in $M_\G$, is precisely the Donaldson--Uhlenbeck
compactification of the moduli space of $\G$--stable bundles. Indeed,
in the case when $\G$ is {\it trivial} this is the result of Li and
Morgan.  \erem

\brem Note that this is not a categorical quotient since $\m M$ is not
ample and is only {\it semi--ample} (Cor \ref{lepot}), i.e some power
of $\m M$ is generated by sections.  \erem

\brem The reduced scheme has a {\it weak categorical quotient}
property for families parametrised by {\it reduced} schemes.\erem

\subsubsection{\it Double duals, associated graded}
Let $F$ be a $\mu$--semistable $\G$--torsion free sheaf over $Y$.  Let
$gr^\mu_\G(F)$ be the graded torsion free polystable sheaf associated
to it's Jordon--Holder filtration. Let $F^{**}$ denote the double dual
of $gr^\mu_\G(F)$; it's a polystable bundle (since $Y$ is a surface, a
reflexive sheaf is locally free). Let $l_F:Y \rightarrow \mathbb N$ be
the function given by $x \mapsto l(F^{**}/gr^\mu(F))_x$, which
associates an element in $S^l_\G(Y)$(length $l$ $\G$--cycle) with
$l=c_2(F)-c_2(F^{**})$. We denote by $Z_F$ the $0$--cycle:
\[
Z_F := \sum_{x \in Y} l(F^{**}/gr^\mu_\G(F))_{x} \cdot x 
\]  
Both $F^{**}$ and $Z_F$ are well defined, i.e. they do not depend on
the choice of filtration.

\subsection{\it Points of the moduli}

The main aim of this subsection is to describe the points of the
moduli space $M_\G$. Towards this we have the following theorem.

Let $Quot(E,l)$ denote the Quot scheme which parametrises all
$0$--dimensional quotients of $E$ of length $l$, where $E$ denotes an
arbitrary torsion--free sheaf on $Y$. If $E$ is a $\G$-vector bundle
on $Y$ the scheme $Quot(E,l)$ gets a natural $\G$--structure and we
can again consider the closed subscheme of $\G$--fixed points in
$Quot(E,l)$. We denote this closed subscheme by $Quot^{\G}(E,l)$.
Clearly this scheme parametrises $0$--dimensional $\G$--quotients of
$E$ of length $l$. 

The $l$--fold symmetric product $S^l(Y)$ parametrises $0$--cycles on
$Y$ of length $l$; again, since $Y$ is a $\G$--surface, by taking the
fixed point subscheme we get the scheme $S^l_\G(Y)$ of zero
dimensional $\G$--invariant cycles of length $l$ on $Y$.  There is
universal sheaf exact sequence on $Y \times Quot(E,l)$:
\beqa\label{tautological}
\xymatrix{
0 \ar[r] & \m E \ar[r] & \m O_{Quot} \otimes E \ar[r] & {\m T} \ar[r] & 0 
}
\eeqa
where $\m E$ is a flat family of torsion--free sheaves on $Y$
parametrised by $Quot(E,l)$. Similarly, we have a $\G$--invariant
exact sequence on $ Y \times Quot^{\G}(E,l)$ with $\m E$ a family of
$\G$-invariant torsion--free sheaves on $Y$.

\beqa
\xymatrix{
Quot^{\G}(E,l) \ar[r]_{\mbox{\tiny {inclusion}}} \ar[d]_{\psi_\G} &
Quot(E,l) \ar[d]_{\psi} \\
S^l_\G(Y) \ar[r]_{\mbox{\tiny {inclusion}}} & S^l(Y) 
}
\eeqa

\brem If $F$ is $\G$--semistable torsion free sheaf we can construct a
family $\m F$ parametrised by $\mathbb P^1$ such that $\m
F_{\infty}=\gF$ and $\m F_t=F$ for all $t \in \mathbb P^1 - \infty$.
This means that $\phi(F) = \phi(\gF)$, where $\phi :\Rss \rightarrow
M_{\G}$ is the canonical morphism. Hence we can restrict to
polystable case alone. It is easy to see that double dual of any
$\G$--sheaf gets a canonical $\G$--structure.\erem


\brem\label{reptypes} Consider the closed subvariety $S^l_\G(Y)$ of
$\G$--invariant cycles on $Y$. Let $Z \in S^l_\G(Y)$ and write $Z =
\sum m_i y_i$. Then the points $y \in Supp(Z)$ can be of the
following types:
\begin{enumerate}\label{torsiontypes}
\item A point $y \in (X \setminus \goth D)$, where ${\goth D}_{Y/X} = {\goth D}$ is
  the ramification divisor of the covering map $p: Y \to X$.
  
\item A general point $y$ contained in an irreducible component
  $(p^*D)_{red}$, the isotropy subgroup $\G_y$ being the cyclic group
  $\G_n$ of order $n$ determined by the irreducible component
  containing $y$.

\item A general point $y$ of an irreducible component of the
  ramification divisor for $p$ not contained in $(p^*D)_{red}$.

\item A special point $y$ contained in $(p^*D)_{red}$, the isotropy
  subgroup $\G_y$ of which contains the cyclic group $\G_n$ of order
  $n$ determined by the irreducible component containing $y$.

\item A special point $y$ of the ramification divisor for $p$ not
  contained in $(p^*D)_{red}$.

\end{enumerate}

Consider a torsion-sheaf $T$ supported at $y \in Supp(Z)$ of length
$m$. Then we can consider the vector space $V$ of its section of
dimension $dim(V) = m$. We view the vector space $V$ endowed with a
${\G}_y$--module structure. For $T_{my}$ to be a quotient of a
$\G$--bundle $E$ on $Y$ of local type $\tau$, the $\G_y$--module
structure on $V$ will have constraints imposed on it arising from the
$\G_y$--module structure on $E|_{U_y}$ which has already been
described in \eqref{localtype}.

Let $Z \in S^l_\G(Y)$ and write $Z = \sum m_i y_i$. For each torsion
sheaf $T_Z$ with support $Z$, fixing a $\G$--structure is equivalent
to fixing a tuple of representations $({\rho}(y_i))$ with $\rho(y_i) :
{\G_{y_i}} \to GL(V)$. Moreover, for any $\gamma \in \G$, since
$\gamma y_i \in Supp(Z)$, we further need that the representation
$\rho(\gamma y_i)$ is the $\gamma$--conjugate to $\rho(y_i)$. \erem

\begin{note}\label{cyclenotation}For a given tuple of representations 
  $\rho(y_i)$ associated to the points in the support of the cycle
  $Z$, we attach a label to the $\G$--cycle $Z$ and denote it by
  $Z(\rho(y_i))$. So an equality $Z_{F_1}(\rho(y_i)) =
  Z_{F_2}(\rho(y_i))$ means that the support of the cycles coincide
  and the torsion sheaves $T_{Z_1} \simeq T_{Z_2}$ are identified as
  $\G$--torsion sheaves. \end{note}

\bth\label{pointsofmoduli} Let $F_i$, $i = 1,2$, be two
$\mu$--semistable $\G$--torsion free sheaves of rank $r$ on $Y$ with
fixed Chern classes $c_1$ and $c_2$. Then $F_1$ and $F_2$ define the
same point in $M_{\G}^{\mu ss}$ if and only if $F_1^{**} \cong_{\G}
F_2^{**}$ and $Z_{F_1}(\rho(y_i)) = Z_{F_2}(\rho(y_i))$.  \eeth

\brem This theorem is proved after the proofs of Proposition
\ref{points1} and Lemma \ref{points2}. \erem

\bprop\label{points1} Let $E$ be a $\G$--polystable vector bundle as
above. Then the {\em connected components} of the fibres of the
morphism $\psi_\G$ are indexed by the representation tuple
$({\rho}(y_i))$ as discussed above in Remark \ref{reptypes}. \eprop

\pf Consider $Z \in S^l_\G(Y)$ and let $T_Z$ be the torsion sheaf with
support $Z$. Let $y \in Supp(Z)$ and lets its multiplicity in $Z$ be
$m$. We first observe that for any $\G$--torsion free sheaf $F \in
{\psi}_G^{-1}(Z)$ canonically induces a tuple of representations
${\rho}(y_i)$ for each of the points $y_i \in Supp(Z)$.

For the given decomposition of $Z$ let us denote a given
representation type on the torsion sheaf $T$ by $T(\rho)$. In other
words, we fix the representation types on $T$ for each point $y \in
Supp(Z)$.

Consider a $\G$--quotient $q: E \to T_Z(\rho)$. We first reduce the
study of such quotients to a local question. 

\begin{itemize}
\item Since $Z$ is a $\G$--cycle, if $y \in Supp(Z)$ so does $\gamma
  y$ for each $\gamma \in \G$. Furthermore, the multiplicities $m$ at
  $y$ and $\gamma y$ also coincide.

\item Giving a $\G$--structure on $T_Z$ is therefore giving
  $\G_y$--structure to $T_{my}$ such that at $\gamma y$, the
  $\G_{\gamma y}$--structure is conjugate to the one at $y$.
  
\item Again, since $E$ is a $\G$--bundle, for any $y \in Supp(Z)$,
  there is a $\G_y$--invariant analytic neighbourhood $U_y$ as in
  \eqref{torsiontypes} such that $E|_{U_y}$ is associated to a
  representation $\G_y \to GL(r)$. Furthermore, at $\gamma y$ for each
  $\gamma \in \G$, the local representation is conjugate to the one at
  $y$ by the element $\gamma$.

\item Giving a $\G$--quotient $q$ as above implies giving quotients
  $q_i : E_i \to T_i(\rho(y_i))$, and where $E_i$ are bundles
  restricted to neighbourhoods of the points in the support of $Z =
  \sum_i m_i y_i$ and $T_i(\rho(y_i)) = T_{m_i y_i}$ with a fixed
  $\G_{y_i}$--module structure on the torsion sheaf $T_{m_i y_i}$.
  Further, the quotient map at $\gamma y_i$ is conjugate to the one at
  $y$.

\item Thus, the problem of studying $\G$--quotients reduces to the
  study of $\G_y$--quotients in a $\G_y$--invariant neighbourhood of
  $y$ . In other words, such a quotient is a point in the product of
  {\it equivariant punctual quot schemes} which we describe below.

\end{itemize}

We therefore need to handle the various points in the possible
singular loci of $\G$--torsion free sheaves as listed in
\eqref{torsiontypes}. 
 
For any point $y \in Supp(Z)$ with multiplicity $m$, suppose that
$\rho(y) : \G_y \to GL(V)$ is already fixed with $dim(V) = m$. Let $V
= \oplus_l~~b_l V(l)$ be the isotypical decomposition as a
$\G_y$--module, with $V(l)$ denoting irreducible $\G_y$--modules.

Consider $E|_{U_y}$ where $U_y$ is an analytic neighbourhood of $y$ as
in \eqref{localtype}. Since the bundle $E|_{U_y}$ is associated to a
representation $\G_y \to GL(r)$, we get an isotypical decomposition
$E|_{U_y} \simeq \oplus_l ({\mathcal O}_{U_y}^{a_l} \otimes V(l))$.

Then, giving a $\G_y$--quotient $q : E|_{U_y} \to T_{my}$ imposes some
natural constraints on $V$, namely, that the $V(l)$'s that occur in
$V$ as a $\G_y$--module must also occur in $E|_{U_y}$ with obvious
bounds on the $a_l$ and $b_l$. With this out of the way, giving $q$ is
equivalent to giving quotients
\[
q_{a_l,b_l} :{\mathcal O}_{U_y}^{a_l} \to T_{b_l y}
\]
twisted by $Id|_{V(l)}$, for each $V(l)$ occurring in $V$.

Since $q_{a_l,b_l}$ is a torsion quotient without any $\G_y$--action,
the irreducibility of the {\it punctual} quot scheme $Quot({\mathcal
  O}_{U_y}^{n_l},m_l)$ is immediate by the results of Jun Li \cite{li}
Baranovsky \cite{bara}, Ellingsrud--Lehn \cite{elli}. Note that we
have this since $Y$ is smooth.

The case when $\G_y$ is trivial, i.e where $y$ avoids the ramification
is easy to handle. In fact, in this case it follows immediately by the
old result quoted above. Thus by the above discussion, it follows that
the {\it equivariant punctual quot scheme} is also {\em irreducible}.

This implies that, fixing the representation type for the torsion
sheaf $T_Z$ gives a connected component of the fibre of $\psi_\G$.

\bcor\label{diffcomp} Let $F_1$ and $F_2$ be two $\G$-polystable
torsion--free sheaves obtained as kernels of two maps in
$Quot^{\G}(E,l)$ and lying in the same fibre of the map $\psi_\G$.  If
we have a $\G$--isomorphism $F_1^{**} \cong_{\G} F_2^{**}$, then $F_1$
and $F_2$ give the same point in the moduli space if and only if they
lie in the same component of the fibre of $\psi_G$ given by a
representation tuple $\rho(y_i)$. \ecor

\pf The fact that $F_i$ ($i = 1,2$) both correspond to points in
$Quot^{\G}(E,l)$, and the assumption that $F_1^{**} \cong_{\G}
F_2^{**}$ implies that we have
\[
E \cong_{\G} F_1^{**} \cong_{\G} F_2^{**}
\]
with the $\G$--structure on $E$ fixed before. 

Let $F_1$ and $F_2$ be (non--uniquely) represented by a two closed
points $q_i \in Quot^\G(E,l)$, $i=1,2$. We think of $F_i$ themselves
as points in $Quot^\G(E,l)$ when there is no confusion.

If $F_1$ and $F_2$ are in a component $S(\rho) \subset
{\psi_\G}^{-1}(Z)$. The line bundle $\m L^N$ is trivial on the fibre
$\psi_\G$ and hence on each component $S(\rho)$ of the fibre of
$\phi_\G$ (since it is the restriction of the determinant bundle on
the fibre of $\psi$). Hence $F_1$ and $F_2$ go to same point in the
moduli space. Conversely, if $F_1$ and $F_2$ lie in different
components, since the line bundle $\m L^N$ is trivial on each
component, one can clearly separate the points $F_i$ by sections of
$\m L^N$. In other words, they go to distinct points of the moduli
space.
\begin{flushright} {\it q.e.d} \end{flushright}

We need to prove the following lemma to complete the proof of the
converse in Theorem \ref{pointsofmoduli}.

\blem\label{points2} Let $F_1$ and $F_2$ are two $\G$-polystable
torsion free sheaves over $Y$. Let $a \gg 0$ and $C \in \mid a{\Theta}
\mid^\G$ is a general $\G$-curve (which exists by the $\G$-Bertini
theorem in the appendix).  Then $F_1\mid_C~\simeq_\G~ F_2\mid_C$ if
and only if $F_1^{**} \simeq_\G F_2^{**}$, where
$F_i^{**}=(gr^\mu_\G(F_i))^{**}, i=1,2$.  \elem

\pf We choose an integer $a$ so large such that restriction of each
summand of $F_1^{**}$ to any general smooth curve $C \in \mid a {\Theta}
\mid^\G$ is $\G$--stable (see Theorem \ref{orbimehtaramanathan}
below). Now we choose one such $C$ in such a way that it avoids finite
set of singular points of $gr^\mu_\G(F_1)$. We note that
$gr^\mu_\G(F_1) \mid_C$ is a polystable bundle over $C$ hence
\[
(gr^\mu_\G F_1) \mid_C \cong gr^\mu_\G(F_1 \mid_C) =
(gr^\mu_\G(F_1)^{**})\mid_C = F_1^{**}\mid_C
\]
The last equality is due to the fact that ``restriction to $C$'' and
``double duals'' commutes with each other.  Now by uniqueness (upto
isomorphism) of Jordan--Holder Filtration of $\G$--semistable bundle
we get $(gr^\mu_\G(F_1))\mid_C \cong_\G F_1^{**}\mid_C$. This shows
that for a general high degree curve $C \in \mid a {\Theta} \mid^\G $, the
bundles $F_1\mid_C$ and $F_2\mid_C$ are S--equivalent if and only if
$F_1^{**}\mid_C \cong_\G F_2^{**}\mid_C$.

\[
\label{can0}
\xymatrix{
0 \ar[r] &  \m O_Y(-C) \ar[r] & \m O_Y \ar[r] & \m O_C \ar[r] & 0
\\
}
\]

Tensoring the above equation with locally free sheaf $\m Hom(F_1^{**},
F_2^{**})$ one gets the following long exact sequence.

\begin{equation*}
\begin{CD}
0 \rightarrow H^{0}_\G(Y,\m Hom_(F_1^{**},F_2^{**})(-C))
  \rightarrow H^{0}_\G(Y,\m Hom(F_1^{**},F_2^{**})) \rightarrow
\\
H^{0}_\G(Y,\m Hom(Y,\m Hom(F_1^{**},F_2^{**})\mid_C) \rightarrow
  H^{1}_\G(Y,\m Hom(F_1^{**},F_2^{**})(-C)) \rightarrow 
\end{CD}
\end{equation*}

We now observe that since we work over fields of characteristic zero
by Remark \ref{gammacohomology}, we have the following inclusions:
\[
H^{i}_{\G}(Y, E) \subset H^{i}(Y, E)
\]

Using this and the usual Serre duality for sheaves on $Y$, we have:

\[
H^{1}_\G(Y,\m Hom(F_1^{**},F_2^{**})(-C))~~\subset~~H^{1}(Y,((\m
Hom(F_1^{**},F_2^{**})^{*}\otimes K_Y)(C))~~=~~0~~
\]

and similarly,

\[
H^{0}_\G(Y,\m Hom(F_1^{**},F_2^{**})(-C))~~\subset~~H^{2}(Y,((\m
Hom(F_1^{**},F_2^{**})^{*}\otimes K_Y)(C))~~=~~0~~
\]

The vanishing follows by Serre vanishing theorem, since $\m Hom
(F_1^{**},F_2^{**})$ is locally free and $C$ is a high degree curve.

Hence we have

\[
H^{0}_\G(Y,\m Hom(F_1^{**},F_2^{**})) \cong H^{0}_\G(Y,\m
Hom(Y,\m Hom(F_1^{**},F_2^{**})\mid_C).
\]
This implies that $F_1^{**} \mid_C \cong_\G F_2^{**}\mid_C$ if and
only if $F_1^{**} \cong_\G F_2^{**}$.

\enpf

\noindent
{\sf Completion of the proof of Theorem \ref{pointsofmoduli}}.  

So if $F_1^{**} \ncong_\G F_2^{**}$ then two points in $R^{\G}$ goes
to two different points in $M^\G$. Now suppose $F_1^{**} \cong_\G
F_2^{**}$, $Z_{F_1}(\rho(y_i)) \neq Z_{F_2}(\rho(y_i))$; By
\eqref{cyclenotation} this means that either $Z_{F_1}\neq Z_{F_2}$ or
that $Z_{F_1} = Z_{F_2} = Z$, but $F_i$ lie in different connected
components of the fibre of $\psi_{Z}$. 

The second case follows from Cor \ref{diffcomp}. If the cycles
themselves are different then we will show that they go to two
different points.  Observe that we have the following diagram:
\[
\xymatrix{
S^l_\G(Y) \ar[r]_a \ar[d]_b & M_\G \ar[d]_{\phi}\\
S^l(Y) \ar[r]_c & M 
}
\]

By \cite{H} that the map $c$ is a {\it closed immersion}.  Since
$S^l_\G(Y)$ is a closed subset of $S^l(Y)$, it follows that $b$ is
also a closed immersion and hence the composite $c \circ b=\phi \circ
a$ is a {\it closed immersion}. So by our assumption $F_1$ and $F_2$
will go to two different points. This completes the proof of the
converse of Theorem \ref{pointsofmoduli}.

\enpf

To realise the construction as a compactification we need to have the
following proposition.  

\bprop\label{stablembed} The moduli space ${M_\G^{\mu s}}(\m Q)$ of
isomorphism classes of $(\G,\mu)$--stable locally free sheaves with
fixed determinant $\m Q$ on $Y$, is {\it embedded} in the moduli space
$M_{\G}$. \eprop

\pf This follows by Lemma \ref{points2} since $F \simeq F^{**}$ for a
stable bundle $F$. The fact that the inclusion is an embedding can be
ensured by choosing $C$ to be of larger degree.  \enpf

\brem Let $M_{\G}^{\mu s}(r,\mathcal Q,c_2)$ denote the moduli space
of $(\G,\mu)$--stable bundles of rank $r$, fixed determinant $\mathcal
Q$ and second Chern class $c_2$. The closure of this moduli space in
$M_{\G}$ gives the desired Donaldson-Uhlenbeck compactification. This
can {\em set theoretically} be described as a stratified space in
terms of $(\G,\mu)$--polystable {\it bundles} with decreasing $c_2$ as
follows:

\beqa {\overline {M_{\G}^{\mu s}(r,\mathcal Q,c_2)(\tau)}} \subset
\coprod_{l\ge 0, \rho} M_{\G}^{\mu-poly}(r, \mathcal
Q,c_2-l)(\tau)\times {S^l_\G}(Y)(\rho) \eeqa

where $M_{\G}^{\mu-poly}(r,\mathcal Q,c_2)(\tau)$ denotes the subset
representing $\G$--polystable locally free sheaves of type $\tau$ and
${S^l_\G}(Y)(\rho)$ consists of zero cycles $Z(\rho(y_i))$ as in
\eqref{cyclenotation}.

\begin{note}\label{modulinotation} We denote by $M^{\pmb{\alpha}}_{k,\bf
    l,\bf r}$ the {\em moduli space of parabolic stable bundles of
    rank $r$} with specified parabolic datum. The tuple
  $(\pmb{\alpha},k,\bf l,\bf r)$ is defined as follows:
\begin{itemize}
\item $ {\pmb{\alpha}} = (\alpha_1,\alpha_2,...,\alpha_l)$,
\item ${\bf l} = (deg(F_1),deg(F_2),...,deg(F_l))$
\item ${\bf r} = (rank(F_1/F_2),rank(F_2/F_3),...,rank(F_l/F_{l+1})$ 
\item $k$ stands for the second Chern class of a vector bundle.  Here
  we follow the notation in \cite{kron} 
\end{itemize}

\end{note}

Recall the correspondence (\eqref{piparab}) between the {\it
  polystable parabolic bundles} on $X$ with given parabolic datum and
$par c_2 = \kappa$ and $(\G,\mu)$--polystable bundles of type $\tau$
on a Kawamata cover $Y$ (see \ref{parabgamma} and \ref{localtype}).
By the description of the above moduli space ${\overline {M_{\G}^{\mu
      s}(r,\mathcal Q,c_2) (\tau)}}$ we get an intrinsic description
the compactification of the moduli space $M^{\pmb{\alpha}}_{k,
  {\mathbf j}, \bf r}(r, \mathcal P, \kappa)$ {\it set-theoretically}
as a stratified space in terms of moduli space of parabolic
$\mu$-polystable bundles with fixed determinant $\mathcal P$ and with
decreasing $\kappa = par c_2$ as follows:

\beqa {\overline {M^{\pmb{\alpha}}_{k, {\mathbf j}, \bf r}(r, \mathcal
    P, \kappa)}} \subset \coprod_{l\ge 0} M^{\pmb{\alpha}-poly}_{k',
      {\mathbf j}', \bf r}(r, \mathcal P, \kappa - l) \times S^l(X).
    \eeqa
    
    where by $M^{\pmb{\alpha}-poly}_{k, {\mathbf j}, \bf r}(r,
    \mathcal P, \kappa)$, we mean the set of isomorphism classes of
    {\it polystable} parabolic bundles with parabolic datum given by
    $(\pmb{\alpha},\bf l,\bf r, \bf j)$, fixed determinant $\mathcal
    P$ and with topological datum given by $k$ and $\kappa$ as
    mentioned above.\erem

\section{Existence of $\G$-stable bundles}

The aim of this section is to prove the existence of $\G$--stable
bundles of {\it rank two} with the assumption of {\it large} $c_2$ or
what is termed {\it asymptotic non-emptiness}. The bound on $c_2$ is
dependent on the polarisation unlike the result of Taubes and
Gieseker. The strategy is to generalise the classical Cayley-Bacharach
property for $\G$--bundles and prove the non-emptiness along the lines
of Schwarzenberger--Serre in the usual surface case.

We remark that, although the moduli space of parabolic sheaves was
constructed on any smooth projective variety (but with the Gieseker
notion of semistability), to the best of our knowledge, the
non-emptiness of these moduli spaces has not been hitherto
established. In this paper we do this over a surface and also show
that at least one component is generically smooth for large values of
$c_2$. {\em As before, we make the following assumptions throughout
  this section}: $Y$ is a smooth projective $\G$--surface which arises
as a ramified Kawamata cover of the smooth projective surface $X$. Let
$p: Y\lr X :=Y/{\G}$ as before denote the covering morphism.

Let $D$ denote the parabolic divisor and $D \, = \, \sum_{i=1}^c D_i$
be the decomposition of the divisor $D$ into its irreducible
components. Since we will be primarily interested in {\em rank two
  bundles}, we have the following weights:
\[
0 \le \alpha_1 < \alpha_2 < 1
\]
where $\alpha_i = \frac{m_i}{N}$. We fix as above a very ample divisor
$\Theta_1$ on $X$ and let $\Theta = p^*(\Theta_1)$.

\bth\label{nonemptiness} The moduli space ${M_\G^{\mu s}}(2,\m Q)$ of
$\G$-stable bundles of rank two and of type $\tau$ and fixed
determinant $\m Q$, on a smooth projective $\G$--surface $Y$ is
nonempty if $c_2(E) \gg 0$ and if $\alpha_2 < \frac{2\cdot
  {\Theta_1}^2}{\sum D_i \cdot {\Theta_1}}$.  Hence, the moduli space
of parabolic bundles on $X$ of rank two with given quasi--parabolic
structure and with ${par c_2}(V) \gg 0$ is non-empty.  \eeth

\brem The parabolic stable bundle that is shown to exist will depend
on the choice of the polarisation $\Theta_1$ on $X$. \erem

\subsection{\it Orbifold Cayley-Bacharach property}

\brem In this section we make the assumption that
$\G$--line bundles that we work with are of type ${\tau}$. \erem

\bdefe\label{good} Let $Y$ be a smooth projective $\G$--surface. Let
$p: Y\lr X$ be a morphism where $X:=Y/{\G}$ arising from the Kawamata
covering lemma. Let ${\goth D}_{Y/X} = {\goth D}$ be the ramification
locus in $X$ and $R$ be a subset of codimension two consisting of
reduced points of length $l$ such that $R\cap \g D =\emptyset$ in $Y$.
Let $Z = p^{*}(R)$. Then we term the cycle $Z$ in $Y$ a {\sf good}
$\G$--cycle.  \edefe

\brem Let $0 \le \beta < \alpha < 1$. Consider $\G$ line bundles $L =
L^{(\alpha)}$, and $M = M^{(\beta)}$ on $Y$ and let $P= M\otimes
L^*\otimes K_Y^{(\alpha - \beta)}$ (see notation in
\eqref{twistedcanonical}).  

By tensoring the standard exact sequence for the ideal sheaf $\m I_Z$
by $P$ we have $ 0\lr \m I_Z\otimes P\lr P\lr \m O_Z\otimes P\lr 0$.
This induces the following exact sequence of $\G$ cohomology groups:

\beqa
\begin{CD}
  0\lr H^0_{\G}(P\otimes \m I_Z)\lr H^0_\G(P)\lr H^0_\G(P\otimes \m O_Z)
  \lr
  \\
  H^1_\G(P\otimes \m I_Z)\lr H^1_\G(P)\lr H^1_\G(P\otimes \m O_Z)=0
\end{CD}
\eeqa

Let $dim H^0_\G(P)=l_1$. Then by choosing a generic $0$--cycle $Z =
p^{*}(R)$ as above such that $l(Z)>l_1$ it is easily seen that we make
sure $H^1_\G( P\otimes \m I_Z)\neq 0$. This implies that there exists
at least one $\G$--torsion free sheaf $E$ on $Y$ which is a non-split
extension of $M\otimes \m I_Z $ by $L$. \erem

\bdefe\label{ocb} {\sf Let $0 \le \beta < \alpha < 1$, and let $L =
  L^{(\alpha)}$ and $M = M^{(\beta)}$ be two $\G$ line bundles of type
  $\tau$ on $Y$} and $Z$ be a {\it good} $\G$--cycle. We say that the
$\G$--triple $(L,M,Z)$, satisfies the {\it Orbifold Cayley Bacharach}
property, (or in short {\it OCB}) if the following holds: for any
section $s\in H^0_\G(M\otimes L^*\otimes K_Y^{(\alpha - \beta)})$ if
the restriction of $s$ to a good $\G$--cycle $Z' \subset Z$ is zero
implies that $s|_Z=0$, where $Z'\subset Z$ is a {\it good }
$\G$--cycle such that $l(Z')=l(Z)-d$, where $d = |\G|$.  \edefe

Let $Z' \subset Z$ be good $\G$-cycles. Consider the exact sequence of
ideal sheaves: 
\[
0\lr \m I_Z\lr \m I_{Z'}\lr \m O_B\lr 0.
\]
Tensor this exact sequence with $M$. By applying the
$Hom_\G(--, L)$--functor to $0\lr M\otimes \m I_Z \lr M\otimes \m
I_{Z'} \lr M \otimes \m O_B\lr 0$ we get a map
\[
\psi_{Z'} :Ext^1_\G(M\otimes \m I_{Z'}, L)\lr Ext^1_\G(M\otimes \m
I_Z,L)
\] 
of $\G$--extensions. 

\blem\label{extocb} Let $(L,M,Z)$ be a $\G$--triple which satisfies
{\it OCB}.  Then we have:
\[
\cup Image(\psi_{Z'})\neq Ext^1_\G(M\otimes \m I_Z,L)
\]
for all good $\G$--cycles $Z'\subset Z$ with $l(Z')=l(Z)-d$. \elem

By tensoring the exact sequence $0\lr \m I_Z\lr \m I_{Z'}\lr \m O_B\lr
0 $ with $P=M\otimes L^*\otimes K_Y^{(\alpha - \beta)}$ we get the
following exact sequence:

\beqa
\begin{CD}
0\lr H^0_\G(P\otimes \m I_Z) \lr H^0_\G(P\otimes \m I_{Z'})\lr
H^0_\G (P\otimes \m O_B) \lr 
\\
H^1_\G(P\otimes \m I_Z)\lr
H^1_\G(P\otimes \m I_{Z'})\lr H^1_\G(P \otimes \m O_B)=0 
\end{CD}
\eeqa 

Here we note that the assumption that the triple $(L,M,Z)$ satisfies
{\it OCB} implies that $H^0_\G(P\otimes \m I_Z)\cong H^0_\G(P\otimes
\m I_{Z'})$. Therefore by dualizing we have: 
\[
0 \lr H^1_\G(P\otimes \m I_{Z'})^*\lr H^1_\G(P\otimes \m I_Z)^* \lr V \lr
0
\]
where $V$ is the complex vector space invariant under $\G$ which is
precisely the dual of the space of sections of the torsion sheaf
$H^{0}_\G(P \otimes \m O_B)$. Note that $V$ is independent of $Z'
\subset Z$ and depends only on $l(Z')$. This in particular implies
that $ H^1_\G(P\otimes \m I_{Z'})^* \subsetneq H^1_\G(P\otimes \m
I_Z)^*$.

Since the finite union of {\it proper} subspaces of finite dimensional
vector spaces is not equal to the vector space (we are over an
infinite field!) we have $ \cup H^1_\G(P\otimes \m I_{Z'})^*\neq
H^1_\G(P\otimes \m I_Z)^*$. The lemma now follows by Serre duality
(Lemma \ref{serreduality}), which gives the identification
$Ext^1_\G(M\otimes \m I_Z,L) \simeq H^1_\G(P\otimes \m I_Z)^*$.
\begin{flushright} {\it q.e.d} \end{flushright}

\blem Let $(L,M,Z)$ be a $\G$--triple which satisfies {\it OCB}. Then
for $l(Z) \gg 0$, there exists a $\G$--extension
\[
0 \lr L \lr E \lr M \otimes \m I_Z \lr 0
\]
with $E$ {\sf locally free}.\elem

\pf Suppose now that $E$ is not locally free. This implies that the
set $Sing(E)$, namely the singular locus of $E$ where $E$ fails to be
locally free, is a $0$-cycle $A\subset Z$, where $A$ is a $\G$--cycle.
Let $a\in A$ then $p^{-1}(p(a))=\sum \gamma .a = B\subset A.$ 

Let $T_A$ denote the torsion sheaf supported at Sing(E). Note
that we have an inclusion of torsion sheaves $T_B\subset T_A.$
Therefore we get the following commutative diagram of $\G$ torsion
free sheaves on $Y$.  \beqa
\begin{CD}
  \xymatrix{ 0 \ar[r] & E \ar[r] \ar@{=}[d] & E^{**} \ar[r] & T_A
    \ar[r] & 0
    \\
    0 \ar[r] & E \ar[r] & E' \ar[u] \ar[r] & T_B \ar[u] \ar[r] & 0\\
} 
\end{CD}
\eeqa where $E' $ be the corresponding subsheaf of the $E^{**} $ to
$\m O_B$.  Note that since $L$ is locally free the saturation of $L$ in
$E'$ is $L$ itself.  

We therefore obtain an extension $E'$ of $M\otimes \m I_{Z'}$ by $L$
using the above commutative diagram where $Z'$ is the $\G$ cycle
corresponding to the {\it good} cycle $R'\subset R$ induced by the set
$A - B$ and $l(Z')= l(Z)-d$ where $d$ is the order of the group $\G$.
Also we have the following commutative diagram of $\G$--sheaves on $Y$
given by two $\G$--sheaves $E$ and $E'$.

\beqa \xymatrix{
  && 0 \ar[d] & 0 \ar[d] \\
  0 \ar[r] & L \ar[r] \ar@{=}[d] & E \ar[r] \ar[d] & M \otimes \m I_Z \ar[r] \ar[d] & 0 \\
  0 \ar[r] & L \ar[r] & E' \ar[r] \ar[d] & M \otimes \m I_{Z'} \ar[r] \ar[d] & 0 \\
  && T_B \ar@{=}[r] \ar[d] & T_B \ar[d]  \\
  &&0&0 } \eeqa

It is clear from the above two diagrams $\psi_{Z'}(E') = E$. By
Lemma \ref{extocb} it follows immediately that there exists locally
free sheaves which can be realised as extensions as desired.
\begin{flushright} {\it q.e.d} \end{flushright}

Now we give the construction of rank two {\it $\G$--stable} vector
bundles as a extension of $M\otimes \m I_Z$ by $\m O_Y$ where $M$ is a
$\G$ line bundle on $Y$.

\brem\label{moregoodness} Let $L$ be a $\G$--line bundle on $Y$ and
let $Z$ be a {\it good} $\G$--cycle. Therefore, $Z = p^{*}(R)$ for a
cycle $R \subset X$ of distinct reduced points away from $\g D$. Under
these conditions we observe the following easy fact:
\[
{p^{\G}_*}(L \otimes \m I_Z) \simeq {p^{\G}_*}(L) \otimes \m I_R
\]
\erem

As before, we fix a very ample divisor ${\Theta}_1$ on $X$ and let ${\Theta} =
p^{*}({\Theta}_1)$ (which is therefore an ample divisor on $Y$). All our
degree computations are with respect to these choices.

\subsubsection{\it Classical Cayley-Bacharach}\label{classicalcayley} 
Let $C$ be a divisor on $X$ with $-2 \Theta_1^2 < C \cdot \Theta_1
\leq 0$. Let $Q = 2 \Theta_1 - C$. Then we have the following well
known result:

\blem\label{usualocb} Let $l \geq h^0(X, Q \otimes K_X)$. Then for a
{\it generic} $0$--cycle $R$ in $Hilb^{l+1}(X)$ we have the usual
Cayley-Bacharach property for the triple $(\m O_X, Q, R)$.\elem

\pf For the sake of completeness we briefly indicate a proof. We first
observe that for generic choice of $T \in Hilb^l(X)$, $l \geq h^0(X, Q
\otimes K_X)$ implies $h^0(X, Q \otimes K_X \otimes {\m I}_T) = 0$.
Let $V_l \subset Hilb^l(X)$ consist of {\it reduced} $0$--cycles and
\[
U_l = \{T \in V_l | h^0(X, Q \otimes K_X \otimes {\m I}_T) = 0 \}
\]
an {\it open dense subset} of $V_l$. Let $\m T$ be the universal
family in $V_{l+1} \times X$, i.e $\m T = \{(T,x) \in V_{l+1} \times X
| x \in Supp(T) \}$ and consider the surjection $f: \m T \to V_l$,
$f(T,x) = T - x$ and the second projection $p : \m T \to V_{l+1}$.
Observe that $p(\m T - f^{-1}(U_l)) \subset V_{l+1}$ is a proper
closed subset. Choose $R \in V_{l+1} - p(\m T - f^{-1}(U_l))$ implying
$p^{-1}(R) \subset f^{-1}(U_l)$ i.e $\forall x \in Supp(R)$, $(R - x)
\in U_l$, hence $h^0(X, Q \otimes K_X \otimes {\m I}_{R - x}) = 0 $,
$\forall x \in Supp(R)$.\enpf

\brem\label{strongCB} In fact, we observe that this choice of $l$
forces something stronger, namely $H^0(Q \otimes K_X \otimes \m I_{R})
= 0$. Moreover, for any $x \in Supp(R)$ we even have $H^0(Q \otimes
K_X \otimes \m I_{R - x}) = 0$ {\it which implies the Cayley-Bacharach
  property}. So if both these vanishings hold, we term the triple $(\m
O_X, Q, R)$ to have the {\sf stronger Cayley-Bacharach property}.
\erem

\blem\label{specialcycle} There exists a {\em good} $\G$--cycle
$Z_1=p^*(R_1)$ in $Y$ with $l(R_1)\ge 4{\Theta}_1^2$ having the
following property:
\begin{center}
  if $\m L$ is any $\G$--line bundle on $Y$ such that $ h^0_\G(\m
  L\otimes \m I_{Z_1}))>0 $ then $c_1(\m L)\cdot {\Theta} \ge 2 {\Theta}^2$.
\end{center}
\elem

\pf Let $C_1$ and $C_1$ be two smooth curves in $|{\Theta}_1|$ in $X$.
Choose a set $S_1$ of $2{\Theta}_1^2$ distinct points in $S_1 \subset (C_1 -
C_2)$ away from $\g D$ the ramification divisor in $X$. Choose
similarly a set $S_2 \subset (C_2 - C_1)$.

Let $R_1 = S_1 \cup S_2$ and let $Z_1 = p^*(R_1)$. Suppose that we
have $h^0_\G(\m L \otimes \m I_{Z_1})>0$. Then from the above Remark
\ref{moregoodness} we get $h^0(p_*^{\G} (\m L)\otimes \m I_{R_1})>0$.

Let $ p^{\G}_*(\m L) = \m L'$. Observe that $\m L$ and $\m L'$ are
both effective. By an abuse of notation, we will continue to denote by
$\m L$ and $\m L'$ divisors in the linear equivalence of the line
bundles.

Suppose that the effective divisor $\m L'$ contains $C_1$ and $C_2$ as
its components. Then
\[
c_1(\m L')\cdot {\Theta}_1 \ge 2{\Theta}_1^2.
\]
If $\m L'$ does not have $C_i$ for some $i=1,2$ then we have 
\[
c_1(\m L') \cdot {\Theta}_1 = \m L'\cap C_i\ge l(S_i) = 2{\Theta}_1^2.
\]
Therefore $c_1(\m L') \cdot {\Theta}_1 \ge 2{\Theta}_1^2$.  Now 
\[
c_1(\m L) \cdot {\Theta} = {deg}_Y(\m L) = (pardeg(p_*^\G(\m L))~|\G|\ge
{deg}_X(\m L')~|\G| \ge 2{{\Theta}_1}^2~|\G| = 2 {\Theta}^2.
\]
\begin{flushright} {\it q.e.d} \end{flushright}

\brem\label{existenceofocb} Let $\m Q \in Pic(Y)$ be a $\G$--line
bundle obtained as follows: Let $Q$ be a line bundle on $X$ and
consider $\m Q \simeq p^*(Q) \otimes \m O_Y^{(\alpha_2)}$, where by
$\m O_Y^{(\alpha_2)}$ we mean the trivial bundle $\m O_Y$ with a
$\G$--structure of type $\tau$ given by multiplication by the
character corresponding to $\alpha_2$ (see \eqref{localtype} for
notation).

Let $0 \le \alpha_1 < \alpha_2 < 1$. Then we claim that for a suitable
choice of $Q$ on $X$, we can ensure that the triple $(\m
O_Y^{(\alpha_1)},\m Q, Z) $ satisfies the orbifold Cayley Bacharach
property with respect to the cycle $Z$. By definition $Z = p^*(R)$. So
we need simply choose $Q$ on $X$ such that the triple $(\m O_X, Q, R)
$ has the usual Cayley-Bacharach property which we get by
\eqref{usualocb}.  This will involve the choice of {\it generic} $R$
with $l(R) \gg 0$ since we need to avoid the ramification locus. We
choose $R$ and $Q$ with the bounds given by Remark
\ref{classicalcayley} which clearly does the job.  \erem

\subsubsection{\it Choice of $\m Q$ and degree bounds}\label{ineqs} 
Let $\gamma = \alpha_2 \cdot \sum deg_X(D_i)$, where $D_i$ are the
irreducible components of the {\it parabolic divisor}.  We let $Q = 2
\Theta_1 - C$, with
\[
-2 \Theta_1^2 + \gamma < C \cdot \Theta_1 \leq 0
\]
{\it This imposes a condition on the weight $\alpha_2$ which we
  therefore have as hypothesis in} Theorem \ref{nonemptiness} (compare
with \eqref{classicalcayley}).

Let $\m Q = p^*(Q) \otimes {\m O_Y^{(\alpha_2)}}$ as in
\eqref{existenceofocb}. Hence

\beqa\label{ineq1}
\frac{c_1(Q) \cdot \Theta_1}{2} < 2\Theta_1^2
\eeqa

Let $d = |\G|$. Then we see that by comparing degrees, we have:
\[
c_1(\m Q) \cdot {\Theta} = {deg}_Y(\m Q) = (pardeg(p_*^\G(\m Q))~d = 
\{{deg}_X(Q) + \gamma \}~d
\]
The non-trivial contribution of $\gamma$ occurs since $p_*^\G(\m Q)$
is a parabolic line bundle with underlying line bundle $Q$ but with
{\it non-trivial parabolic structure}.

Again, since $deg_X(Q) = 2 \Theta_1^2 - C \cdot \Theta_1$, by the bounds
for $C \cdot \Theta$ fixed above and an easy computation gives:

\beqa\label{ineq2}
\frac{c_1(\m Q) \cdot {\Theta}}{2} < 2{\Theta}^2 
\eeqa

\blem Let $\m Q\in Pic(Y)$ a $\G$--line bundle of type ${\tau}$ as in
\eqref{existenceofocb} and \eqref{ineqs} with $\alpha_2$ as in Theorem
\ref{nonemptiness}. Then there is a $\G$--stable rank two vector
bundle $E$ of type $\tau$ with weights $(\alpha_1, \alpha_2)$, with
$det(E)\cong \m Q$ and $c_2(E)=c$.  \elem

\pf First we start with a a triple $(\m O_Y^{(\alpha_1)},\m Q, Z_2)$
which satisfies the orbifold Cayley-Bacharach property. This exist by
what we have already seen (by \eqref{existenceofocb} and
\eqref{classicalcayley}). We in fact choose a $0$--cycle $R_2$ in $X$
to satisfy the stronger property as in \eqref{classicalcayley} and
\eqref{strongCB} and let $Z_2 = p^*(R_2)$.

This gives us a $\G$ locally free extension
$E'$ of $Q\otimes \m I_{Z_2}$ by $\m O_Y$.

Now we choose a {\it good} $\G$--cycle $Z_1$ as in Lemma
\ref{specialcycle} and let 
\[
Z = Z_1 \cup Z_2.
\]  

Then we observe that the triple $(\m O_Y,\m Q, Z) $ also satisfies a
orbifold Cayley Bacharach property. This can be seen as follows: if $Z
= p^*(R)$, then by \eqref{existenceofocb}, its enough to see that $(\m
O_X, Q, R)$ has the usual Cayley-Bacharach property. This immediate,
for if $x \in Supp(R) = Supp(R_1) \cup Supp(R_2)$, then its easy to
see that $H^0(Q \otimes K_X \otimes \m I_{R - x}) = 0$ since we have
assumed the stronger Cayley-Bacharach property for $R_2$ and moreover,
$\m I_{R - x} \subset \m I_{R_2 - x}$ or $\m I_{R - x} \subset \m
I_{R_2}$ depending on whether $x \in Supp(R_2)$ or not.

Therefore we get a new $\G$--locally free extension $E$:
\[
0 \lr \m O_Y^{(\alpha_1)} \lr E \lr {\m Q}^{(\alpha_2)}\otimes \m I_Z \lr 0
\]
We now {\it claim} that any such $E$ is $\G$--stable.  

To see this, consider any $\G$--line subbundle $L$ of $E$. If $L$ is
non-trivial, then composing the inclusion $L \hra E$ with the map
$E\lr \m Q\otimes \m I_Z$ we get a nontrivial $\G$--map $f: L\lr \m
Q\otimes \m I_Z$.  This gives a non-zero $\G$--section
\[
s \in H^0_\G(\m Q\otimes L^*\otimes \m I_Z).
\]  
In particular, $h^0_\G(\m Q\otimes L^*\otimes \m I_{Z}) > 0$ and as a
result $h^0_\G(\m Q\otimes L^*\otimes \m I_{Z_1}) > 0$.  Therefore by
Lemma \ref{specialcycle} we conclude that 
\[
(c_1(\m Q) - c_1(L))\cdot {\Theta}\ge 2 {\Theta}^2.
\]
Hence, $\mu(L) = c_1(L) \cdot {\Theta} \le (c_1(\m Q)\cdot
{\Theta}-2{\Theta}^2)$.  But we know that $\mu(E)= \frac{(c_1(\m
  Q)\cdot {\Theta}}{2}$.  By \eqref{ineq2} we thus have:

\[
\mu(L) \le c_1(\m Q)\cdot {\Theta}-2{\Theta}^2 <\frac{(c_1(\m Q)\cdot
  {\Theta}}{2} = \mu(E).
\]

Hence $E$ is $\G$--stable and clearly of determinant $\m Q$.

Regarding the type of the $\G$--stable bundle $E$ of rank two
constructed above, we observe that we work with a zero cycle $Z$
coming from the complement of ramification divisor. So the action of
$\G$ on $Z$ is a free action. So it does not affect the {\sf type} of
the extension bundle we constructed. 

Therefore, since we start with $\G$--line bundles of type $\tau$ (see
\eqref{localtype}), by giving a type $\tau$ structure to $\m
O_Y^{(\alpha_1)}$ i.e the trivial bundle $\m O_Y$ with the action of
generic isotropies along the irreducible components of the divisor by
the character $\alpha_1$ and similarly $\m Q = p^{*}(Q) \otimes \m
O_Y^{(\alpha_2)}$. Then we get a rank two stable $\G$--vector bundle
of type $\tau$ via the extension:
\[
0 \lr \m O_Y^{(\alpha_1)} \lr E \lr {\m Q}\otimes \m I_Z \lr 0
\]

\begin{flushright} {\it q.e.d} \end{flushright}

\bcor\label{smoothpoints} There exists $(\G,\mu)$--stable on $Y$ with
vanishing obstruction space. \ecor

\pf To see this we make a few easy observations:
\begin{enumerate}
\item The obstruction space of a $\G$--bundle on $Y$ can be easily
  seen to be the space $Ext^{2}_{\G}(E,E)_0$, where the subscript
  stands for the {\it trace zero} part.
  
\item Now we compute the $Ext_{\G}$ using the construction of $E$ as a
  $\G$--extension. The argument is exactly as in \cite[Remark
  5.1.4]{H} and we only use the Remark \ref{gammacohomology} to get
  the vanishing when we make degree and length large.
\end{enumerate}

\enpf










\section {Kronheimer-Mrowka results revisited}
In the section, for the sake of simplicity, we work with $D \subset X$
an irreducible divisor as the parabolic divisor. The other notations
are as in Section 2.

\subsubsection{\it Calculation of the second parabolic Chern class}\label{action}
\blem Consider a general parabolic bundle $(E_*,F_*,\alpha_*)$. Then
we can compute the parabolic Chern classes $E_*$ using the following
formula on $X$. Let us assume that $deg(F_i)=l_i$ with corresponding
weights $\alpha_i$ and $r_i = rank(F_i/F_{i+1})$. Then
\[
par c_1(E) = c_1(E) + (\sum_{i=1}^{i=l} r_i\alpha_i) D 
\] 
and

$par c_2(E) =$
\[ 
c_2(E) + \sum_{i=1}^{i=l}r_i\alpha_i (c_1(E).D) -
\sum_{i=1}^{i=l} \alpha_i(l_i-l_{i+1}) + {\frac
  {1}{2}}\{(\sum_{i=1}^{i=l}r_i \alpha_i) \cdot (\sum_{j=1}^{j=l}r_j
\alpha_j) - (\sum_{i=1}^ {i=l}r_i \alpha_i^2) \}D^2
\]
\elem 

\pf We can assume without loss of generality that $E$ is a parabolic
direct sum of line bundles $(L_i,\alpha_i)$. It is easy to see that
$F_i/F_{i+1}= \oplus_{j\in J} L_j|D$ with $\alpha_j=\alpha_i$ and
$J\subset I$ where $E=\oplus_{i\in I}L_i$. Then $parc_2(E)=\sum_{i<
  j}(c_1(L_i)+\alpha_i D) (c_1(L_j)+\alpha_j D)$.

Hence

$parc_2(E)=$
\[
\sum_{i<j}c_1(L_i)c_1(L_j)+\sum_{i\neq j}c_1(L_i)\alpha_jD+
\sum_{i<j}\alpha_i\alpha_jD^2
\]
The first term in the above equation is $c_2(E)$.  In the above
equation $\alpha_i$ repeated $r_i$ times and $\sum r_i=r$. We write

$\sum_{i\neq j}^rc_1(L_i)\alpha_jD =$
\[
\sum_{i=1}^r\alpha_i\sum_{j=1}^rc_1(L_j)D-\sum_{i=1}^r\alpha_ic_1(L_i)D = 
\sum_{i=1}^lr_i\alpha_ic_1(E)D-\sum_{i=1}^l\alpha_ic_1(F_i/F_{i+1})
\]
where $l$ is the length of the filtration. So, we get the required
second term of the formula. For the third term we just note that
$\sum_{i\neq j}\alpha_i\alpha_j=2\sum_{i<j}\alpha_i\alpha_j$ and by
usual manipulation we get the above formula. \enpf

As in \cite{kron}, we work with a parabolic vector bundle $E$ of rank
two on $(X,D)$ where $D$ is an irreducible smooth divisor with $c_1(E)
= 0$ and a filtration $0\subset \mathcal L \subset E|_{D}$ with a
single weight $\alpha$ associated with a line subbundle $\mathcal L$.
When $E = L\oplus L^*$ with $c_1(E)=0$ and a filtration $0\subset
\mathcal L \subset E|_{D}$ we get $par c_1(E)=0$ and
\[
par c_2(E) = c_2(E)+ 2\alpha \cdot l - \alpha^2 D^2
\]
where $(- l) $ is the degree of the line bundle ${\mathcal L}$
and $\alpha$ is a corresponding weight.

\subsubsection{\it The boundary points and action}

\cite[Theorem 8.21]{kron} says that there is a one-to-one
correspondence between the set of irreducible connections in the
moduli space $M_{k,l}^\alpha (X,D) $ of $\alpha$ twisted
connections , anti-self dual with respect to the cone-like metric
determined by $\omega$, with holonomy parameter $\alpha =a/v$; and the
set of stable parabolic $SL(2,\mathbb C)$ bundles $(\mathcal E,
\mathcal L,\alpha)$ on $X$, with the same weight $\alpha$, satisfying
$c_2(\mathcal E)=k$ and $c_1(\mathcal L)= -l$.

We consider \cite[Proposition 7.1]{kron1}, which is the parabolic
analogue of the Uhlenbeck compactness lemma. This says that if $A_n$
be a sequence of twisted connections in the extended moduli space
$\bar{M}_{k,l}$ over (X,D), and suppose that the holonomy parameters
$\alpha_n$ for these connections converge to $\alpha\in (0,1/2)$. Then
there exists a sub-sequence, which we continue to call $A_n$, and
gauge transformations $g_n \in \mathcal G$ such that the connections
$g_n(A_n)$ converge away from a finite set of points ${x_i}\subset X$,
to a connection $A$. The solution $A$ extends across the finite set
and defines a point in a moduli space $M_{k',l'}^\alpha$.

In \cite{kron} the difference between $(k,l)$ and $(k',l')$ is
accounted for by what {\it bubbles off} at the points where
convergence fails. Thus, for each point of concentration $x_j$ in $X
\setminus D$ there is an associated positive integer $k_j$, and for
points of concentration $x_i$ in $D$ there is an associated pair
$(k_i,l_i)$ so that $k' = k - \Sigma k_i - \Sigma k_j$ and $l' = l -
\Sigma l_i$.  In \cite{kron} it is remarked that there is no complete
interpretation or description of the possible values of the pairs
$(k_i,l_i)$ in the {\it bubbling off}. The key observation made in
\cite{kron} is that the {\it action} $\kappa$ is precisely the
quantity which is seen to decrease in the {\it bubbling off}.

We wish to interpret this phenomenon in the light of the semistable
reduction theorem (see Appendix Theorem \ref{langton}) as well as the
description of the points in the boundary of the Donaldson-Uhlenbeck
compactification constructed in this paper.  

The analogue of the Uhlenbeck compactness lemma in our setting is the
interpretation of the Langton extension in terms of the points of the
boundary, i.e the limit point of the family $E_{(A - p)}$ of parabolic
$\G$ stable sheaves on $(Spec(A) - p)$ of parabolic Chern class $par
c_2$ coming from Langton criterion is identified with a pair
$(E_p,Z_p)$ where the parabolic Chern class of $E_p$ is $par c_2 - s$
where $s$ is the length of the zero cycle $Z_p$. In other words, the
phenomenon of {\it bubbling off} is seen in the decreasing of the
second parabolic Chern class which is precisely the expected
description seen in the light of Donaldson's theorem in the {\it
  non-parabolic} setting. In the case of rank $2$ as in \cite{kron},
what is termed {\it action} and denoted by $\kappa$ is precisely the
second parabolic Chern class. We may therefore interpret the second
parabolic Chern class as the {\it action} in all ranks as seen from
our construction of the Donaldson-Uhlenbeck compactification.


The invariant $par c_2$ captures all the information about the
invariants $(k,l)$ and $(k',l')$ in the notation of \cite{kron}, and
also the relation between them. Indeed, $parc_2$ can be written in
terms of these $k$ and $l$'s as we have seen above.  And since we use
$\G$ bundles on $Y$, we observe that $parc_2$ is able to recover the
information about these numbers as we have described earlier.  The
term {\it action}, denoted by $\kappa$ in \cite{kron} is nothing but
our $parc_2$. Kronheimer and Mrowka define $\kappa_i = k_i + 2\alpha
l_i$, as the action lost at the point of concentration $x_i\in D$.
They also give the relation between $\kappa$ and $\kappa'$ i.e
$\kappa' = \kappa - \Sigma \kappa_i - \Sigma k_j$ , where $\kappa'$ is
the $parc_2$ of the limiting point in our compactification. Here $k_j$
are the instanton numbers associated with the points of concentration
away from $D$.
 
\subsubsection{\it Concluding remarks}
In the sequel to this work (\cite{balasuman}) we prove the {\it
  asymptotic irreducibility, asymptotic normality and generic
  smoothness} of the moduli space of stable parabolic bundles. These
generalise the work of O'Grady and Gieseker-Li for the usual moduli
spaces of stable bundles on algebraic surfaces.

\section{Appendix}

\subsection{\it The Mehta-Ramanathan restriction theorem for orbifold bundles}

The aim of this section is to prove the Mehta-Ramanathan restriction
theorem for $\G$--sheaves. This in particular gives a different proof
of the restriction theorem for parabolic bundles (proven in
\cite{usha}) but for the type of parabolic bundles which arise as
invariant direct images of orbifold bundles. We remark that for the
purposes of the geometric study of the moduli spaces of parabolic
bundles, our results suffice by the {\it yoga} of variation of
parabolic weights.

\subsubsection{\it Remark on $\G$--Bertini}
One has the following version of $\G$--Bertini theorem needed in the
restriction theorem. We omit the proof which is a straightforward
generalisation of the usual case.

\bth($\G$--Bertini) Let $X= Y/\G$. Let us assume that $X$ is smooth
and $\Theta$ is a pull-back of a very ample divisor $\Theta_1$ on $X$. 

Let the closed embedding $Y \subset \mathbb P^n$ be induced by
$\Theta$ i.e $\mathbb P^n$, the projective space determined by
$\mid \Theta \mid$. Then there exists a $\G$ hyperplane $Z \subseteq
\mathbb P^n$, not containing $Y$, and such that the scheme $Z \cap
Y$ is regular at every point.  Furthermore, the set of hyperplanes
with this property forms an open dense subset of $\mid \Theta
\mid^{\G}$.  \eeth


 

\subsubsection{\it The restriction theorem for orbifold bundles}

We have the following $\G$--Mehta Ramanathan restriction theorem from
which the parabolic version follows easily.  

\bth\label{orbimehtaramanathan} ($\G$--Mehta-Ramanathan theorem) Let
$E$ be a $(\G,\mu)$--semistable (resp stable) $\G$--torsion free sheaf
on a smooth projective $\G$--variety such that $X = Y/{\G}$ is also
smooth and projective. Then the restriction $E|_{C_k}$ to a general
complete intersection $\G$--curve $C_k$ of large degree (with respect
to the pull-back line bundle $\Theta$ as in Bertini above) is
$(\G,\mu)$--semistable (resp stable).  \eeth

\pf Since $(\G,\mu)$--semistability for $\G$--sheaves is equivalent
to the semistability of the underlying sheaf, the non-trivial case is
that of stability. The proof can be seen in the following steps:
\begin{enumerate}
\item\label{poply} Let $E$ be $(\G,\mu)$--polystable. Then the
  underlying bundle $E$ is $\mu$--polystable.  In particular, if $E$
  is $(\G,\mu)$--stable the underlying bundle is $\mu$--polystable
  (not necessarily stable).  For, if we start with a $\G$ stable
  bundle $E$ we can construct a socle $F$ of $E$ with $\mu(F) =
  \mu(E)$ which is invariant under all the automorphisms of $E$, in
  particular invariant under the group $\G$. This contradicts the $\G$
  stability of $E$.

\item\label{poply1} By the effective restriction theorem of Bogomolov
  (cf. \cite{H}), for {\it every} complete intersection curve $C$ in
  the linear system $|m \Theta|$ (the number $m$ effectively
  determined), the restriction $E|_{C}$ is polystable.
  
\item\label{poply2} By the $\G$--Bertini theorem, there always exists
  a $\G$--curve in $|m \Theta|$. Thus, the restriction $E|_{C}$ to any
  $\G$--curve is a $\G$--bundle and also $\mu$--polystable. This
  implies that $E|_{C}$ is a $(\G,\mu)$--polystable bundle on $C$.
  For, we take $\G$-socle $F$ of $E|_{C}$ which is again the socle of
  $E|_{C}$.  Now this is $\mu$--polystable proving that $E|_{C} = F$.

\item Observe that if $E$ is $(\G,\mu)$--stable then it is
  $\G$--simple.  Here we note that we are not saying that it is {\it
    simple}. If not, choose a nontrivial $\G$ endomorphism which
  induces a nontrivial $\G$ subbundle of $E$ with $\mu(F)\ge \mu(E)$
  contradicting the $(\G,\mu)$--stability of $E$.
  
\item By the orbifold version of Enriques-Severi it follows that for
  sufficiently high degree $C$ which is also a $\G$--curve, $E|_{C}$
  is also $\G$--simple.
  
\item Hence, if $E$ is $(\G,\mu)$--stable, then for high degree
  $\G$--curve $C$, the restriction is $\G$--simple and
  $(\G,\mu)$--polystable (by \eqref{poply}, \eqref{poply1}, and
  \eqref{poply2} above), and hence $\G$--stable.
\end{enumerate}

\begin{flushright} {\it q.e.d} \end{flushright}

\subsection{Valuative criterion for semistable orbifold sheaves}

Let $S$ be an algebraic variety over $\mathrm k$. We say a
$\G$-coherent sheaf $E$ on $X\times S$($S$ with trivial $\G$ action)
is a {\it family of torsion-free sheaves} on $X$ over $S$ if, $E$ is
flat over $S$ such that for each $s\in S$ the induced sheaf $E_*$ on
${p}^{-1}(s)$ is $\G$-torsion free sheaf on $X$. We say two such
families $E$ and $E'$ are equivalent if there is $\G$ invertible sheaf
$L$ on $S$ such that $E\cong E'\otimes\mathcal {\mathrm p}_{2}^*(L)$.

Our field $\mk$ is algebraically closed. Let $\mk \subseteq R$ be a
discrete valuation ring with maximal ideal $\mathrm m$ generated by a
uniformizing parameter $\pi$ . Let $\mathrm K$ be the field of
fractions of $R$. Consider the scheme $X_R = X\times Spec R $.  Denote
by $X_{\mathrm K} $ the generic fiber and by $X_{\mk}$ the closed
fiber of $X_R$. Let $i$ be the open immersion $X_{\mathrm
  K}\hookrightarrow X_R$ and $j$ be the closed immersion $X_{\mk}
\hookrightarrow X_R$. 

We can now state the main theorem in this section, namely the
semistable reduction theorem for $(\G,\mu)$--semistable torsion-free
sheaves.

\bth\label{langton} Let $E_\mathrm K$ be a $\G$-torsion free sheaf on
$X_\mathrm K$.  Then there exists $\G$--torsion free sheaf $E_\mathrm
R$ on $X_\mathrm R$ such that over $X_\mathrm K$ we have $i^{*}
E_\mathrm R \simeq E_\mathrm K$ and over the closed fibre $X_{\mk}$
the restriction $j^{*}(E_\mathrm R)$ is $(\G,\mu)$--semistable.\eeth

\pf We remark that we need essentially two additional ingredients in
the old proof of Langton to complete our argument. The first one is
that without the demand of semistability by Prop 6 in \cite{langton}
one firstly obtains a {\bf canonical} extension of $E_\mathrm K$ to a
torsion-free sheaf $\tilde {E}$ on $X_\mathrm R$. Since the family
$E_K$ on $X_K$ is given to be a $\G$--sheaf and since the extension is
canonical it follows without much difficulty that the extension also
carries an extended $\G$--action. In other words, the restriction 
$j^{*}(\tilde E)$ to the closed fibre $X_{\mk}$ is also a $\G$--torsion free
sheaf but which could be $\mu$--unstable.

The second step in Langton's proof is to modify the family
successively by carrying out elementary modifications using the first
term of the Harder-Narasimhan filtration (the so-called
$\beta$--subbundle) of the restriction $j^{*}{\tilde E}$. We again
observe that the $\beta$--subbundle being canonical is also a
$\G$--sheaf. In other words, the family remains a $\G$--family even
after the elementary modifications. That the process ends after a {\it
  finite} number of steps is one of the key points in Langton's proof
and we see that we achieve a $(\G,\mu)$--semistable reduction in the
process.

\enpf

\bcor If the generic member of the family $E_K$ is given to be of type
$\tau$ as a family of $\G$--sheaves then so is the closed fibre. \ecor

\pf This is easy to see since the type of the family remains constant
in continuous families. \enpf


\begin{thebibliography}{99}



  
  
\bibitem{balasuman} V. Balaji, S. Bandhopadhyay, Parabolic bundles on
  algebraic surfaces-II, Irreducibility (in preparation)
  
\bibitem{Bi3} V. Balaji, Indranil Biswas, D. S. Nagaraj, On the
  principal bundles over projective manifolds with parabolic structure
  over a divisor, Tohoku Math. J.(2){\bf 53},(2001),no.3,337--367

  
\bibitem{bara} V. Baranovsky, Moduli of Sheaves on surfaces and action
  of oscillator algebra, J.Diff. Geom,{\bf 55},No. 2,(2000), 193-227.
  
\bibitem{usha} Usha Bhonsle, Parabolic sheaves on higher dimensional
  varieties, Math Ann. {\bf 293} (1992), 177-192.
  
\bibitem{biquard} O. Biquard, Prolongement d'un fibr\'{e} holomorphe
  hermitien \'{a} courbure $\mathcal L^p$ sur une courbe ouverte,
  Internat. J. Math. {\bf 3} (1992), 441-453.

  
\bibitem{Bi1} I. Biswas, Parabolic Bundles as Orbifold Bundles, Duke
  Mathematical Journal, {\bf 88} No:2 (1997).

  
\bibitem{Biswaschern} I. Biswas, Chern classes for parabolic bundles,
  J. Math. Kyoto Univ.  {\bf 37} (1997), 597-613
  
\bibitem{Bi4} I. Biswas, On the cohomology of parabolic line bundles,
  Math. Res. Let. {\bf{2}} (1995),783-790.

  
\bibitem{elli} G. Ellingsrud and M. Lehn: Irreducibility of the
  punctual quotient scheme of a surface.  Ark.Mat. {\bf 37},(1999),
  No. 2, 245-254
  
\bibitem{grothendieck} A. Grothendieck, Sur quelques points
  d'alg\`ebre homologiques, Tohoku Math.J.(2),{\bf 9} (1957),
  119--221.
  
\bibitem{groth} A.Grothendieck, Sur la m\'emoire de Weil
  ``Generalisation des fonctions ab\'eliennes'', S\'eminaire Bourbaki,
  Expos\'e 141, (1956-57).
  

  
\bibitem{heinzner} P.Heinzner and F.Kutzschebauch, An equivariant
  version of Grauert's Oka principle, Invent. Math, {\bf 119}, (1995)
  317-346.
  
\bibitem{H} D. Huybrechts and M. Lehn, The Geometry of Moduli Space of
  sheaves, 269 Pages, Friedrick Vieweg \& Son, (1997)

\bibitem{K} Y. Kawamata, Characterization of the abelian varieties ,
  Compositio Math {\bf 43} (1981), 253-276.
  
\bibitem{KMM} Y. Kawamata, K. Matsuda, K. Matsuki, Introduction to the
  minimal model problem in {\it Algebraic Geometry, Sendai,} (1985),
  Adv. Stud. Pure.  Math {\bf 10} North-Holland, Amsterdam, 1987,
  283-360.
  
\bibitem{kron} P. B. Kronheimer and T. S. Mrowka: Gauge theory for
  embedded surfaces, I Topology, {\bf 32},(1993) No. 4, 773-826.
  
\bibitem{kron1} P. B. Kronheimer and T. S. Mrowka: Gauge theory for
  embedded surfaces, II Topology, {\bf 34,}(1995) No. 1, 37-97.


\bibitem{langton} S.Langton: Valuative criterion for families of vector
bundles on algebraic varieties, {\it Annals of Mathematics}(2) {\bf
101} (1975) pp 88-110.

  
\bibitem{Le} J.Le Potier, Fibre determinant et courbes de saut sur
  les surfaces algebriques, Complex Projective Geometry, London
  Mathematical Society: Bergen (1989) 213--240



  
\bibitem{Li} J.Li, Picard groups of the moduli spaces of vector
  bundles over algebraic surfaces, Proc. Symposium. Taniguchi Kyoto
  1994; Lecture Notes in Pure and Applied Mathematics 179 (1996),
  129--146

\bibitem{li} J.Li: Algebraic geometric interpretation of Donaldson's
polynomial invariants of algebraic surfaces, {\it J.Diff.Geom} {\bf
37} (1993), 416-466.
  
\bibitem{linarasimhan} Jiayu Li and M.S. Narasimhan,
  Hermitian-Einstein metrics on parabolic stable bundles.  Acta Math.
  Sin. (Engl. Ser.)  {\bf 15} (1999), no. 1, 93--114.






\bibitem{MY} M. Maruyama, K. Yokogawa, Moduli of parabolic stable
  sheaves, Math. Ann. {\bf 293} (1992),77-99.
  
\bibitem{MS} V.B. Mehta and C.S. Seshadri, Moduli of vector bundles on
  curves with parabolic structure, Math. Ann.  {\bf 248} (1980),
  205--239.


\bibitem{MR1} V.B. Mehta and A. Ramanathan, Restriction of stable
  sheaves and representations of the fundamental group, {\bf 77}
  (1984), 163--172.
  
\bibitem{MR2} V.B. Mehta and A. Ramanathan, Semistable sheaves on
  projective varieties and the restriction to curves, Math. Ann.
  258(1982), 213--226
  
\bibitem{NR} M.S. Narasimhan  and A. Ramanathan, : Openness of the
  semistability condition (unpublished manuscript)
  
\bibitem{S1} C.S. Seshadri, Moduli of $\pi$--bundles over an algebraic
  curve, in Questions on Algebraic Varieties, C.I.M.E, III, Ciclo,
  Varenna (1970), pp 139-260.
  
\bibitem{steer} B. Steer and A. Wren, The Donaldson-Hitchin-Kobayashi
  correspondence for parabolic bundles over orbifold surfaces.  Canad.
  J. Math.  {\bf 53} (2001), no. 6, 1309--1339.


x
  
\bibitem{Bi3Yo} K. Yokogawa, Infinitesimal deformations of parabolic
  Higgs sheaves, Internat. J. Math.{\bf 6}(1995),125-148


\end{thebibliography}
\end{document}